\documentclass[a4paper,12pt]{article}

\usepackage[T1]{fontenc}
\usepackage[utf8]{inputenc}
\usepackage{lmodern}
\usepackage[margin=24mm]{geometry}
\usepackage{amsmath,amssymb,booktabs,tabularx,microtype,xurl,placeins,longtable,graphicx}
\usepackage{tikz}
\usepackage{placeins}
\usepackage{ulem}

\def\fct#1{\mathop{\rm #1}}	

\def\beq{\begin{equation}} 
	\def\eeq{\end{equation}} 
\def\lbeq#1{\begin{equation} \label{#1}} 
	                           
	\def\argmaxx{\fct{argmax}}
	\def\argminx{\fct{argmin}}

\usetikzlibrary{arrows.meta,positioning,shapes.geometric}
\graphicspath{{figures/}}
\usepackage[colorlinks=true,linkcolor=blue,citecolor=blue,urlcolor=blue]{hyperref}
\newcommand{\E}{\mathbb E}

\newcommand{\R}{\mathbb R}
\begin{document}

		\begin{center}
			
			{\Large \bf A particle swarm optimization with a reduced-subspace trust-region safeguard for contact-rich controller tuning} 
			
			\vspace{0.5cm}

			{\large \bf Morteza Kimiaei$^{a}$, Roya Khalili-Amirabadi$^{b}$, 
				Mohsen Jalaeian-Farimani$^{c}$}

			{$^a$ \sl Faculty of Mathematics, University of Vienna, 1090 Vienna, Austria}
			\centerline{\sl Email: \href{mailto:kimiaeim83@univie.ac.at}{kimiaeim83@univie.ac.at}}
			\centerline{\sl WWW: \url{http://www.mat.univie.ac.at/~kimiaei}}\\
			\vspace{0.2cm}
			
			{$^b$ \sl Department of Applied Mathematics, Ferdowsi University of Mashhad, Mashhad, Iran}
			\centerline{\sl Email:
				\href{mailto:roya.khalili.a@gmail.com}{roya.khalili.a@gmail.com}}\\
			\vspace{0.2cm}
			{$^c$ \sl Department of Electronics, Information and Bioengineering (DEIB), Politecnico di Milano, Milan, Italy}
			\centerline{\sl Email:
				\href{mailto:mohsen.jalaeian@polimi.it}{mohsen.jalaeian@polimi.it}}
			\centerline{\sl WWW: \url{https://www.deib.polimi.it/eng/people/details/2012404}}
			\\

		\end{center}
		
		{\bf Abstract.}
	Contact-rich robotic insertion is challenging because
	uncertainty in alignment, friction, and contact conditions
	can lead to unsuccessful assembly, excessive contact forces,
	and jamming. Designing a reliable controller therefore
	requires balancing insertion performance, robustness,
	and recovery from unexpected contact. In this paper,
	we propose \texttt{PSO+TR}, a derivative-free optimization (DFO) method that combines
	the exploratory search of particle swarm optimization (\texttt{PSO}) with a trust-region ({\tt TR}) procedure for local refinement to improve the search process when progress
	stagnates. We apply the method to derivative-free policy optimization of an interpretable
	robotic insertion controller, jointly tuning its motion,
	feedback, and recovery parameters while preserving its
	predefined control structure.
	
The proposed approach is evaluated on 93 eligible benchmark problems from the bound-constrained programs (\texttt{BCP}) collection and a synthetic contact-rich planar
insertion task under uncertain operating conditions.
Three random seeds per benchmark problem yield 279 matched
problem--seed instances per method.
\texttt{PSO+TR} reaches the prescribed target on 228 of these
instances (81.72\%), the highest observed target-attainment rate
among the five compared methods, which also include \texttt{PSO},
randomized black-box optimization (\texttt{VRBBO}), matrix adaptation
evolution strategy (\texttt{MAES}), and differential evolution
(\texttt{DE}).
Across eight robotic optimizer seeds, the deployed controller-selection
pipeline achieves 97.78\% completion under development conditions
and 58.00\% under withheld stress conditions, compared with
91.22\% and 31.12\%, respectively, for the fixed staged controller.
The stress violation rate is reduced from 20.69\% to 2.19\%.
	
	These findings demonstrate the potential of combining
	global and local search for DFO and structured
	contact-rich controller tuning. The robotic results
	are limited to a synthetic planar simulator, and
	further experiments are required to establish
	performance on physical robots.

	

\begin{sloppypar}

\section{Introduction}

Robotic assembly requires motion to remain effective as a robot
interacts with its environment. During insertion, small differences
in alignment, friction, or clearance can determine whether a part
moves freely, contacts a rim, or becomes stuck. Accurate motion
under nominal conditions is therefore only one part of the control
problem. The controller must also accommodate uncertain interaction
conditions while balancing progress and contact loading.

Even when a suitable feedback structure is available, selecting its
numerical parameters can be difficult. Motion amplitudes, controller
gains, and timing decisions jointly influence the resulting behavior.
A setting that shortens completion time for a well-aligned part may
increase contact loading when alignment is poor. Conversely, a more
conservative setting may reduce loading but slow progress or prevent
completion. These interactions make parameter selection a coupled
problem: improving one aspect of the task can degrade another.

The information needed for tuning is also costly to obtain.
Completion, elapsed time, and contact loads emerge over an entire
closed-loop trajectory rather than from a single control update.
Simulation or experimental trials can reveal these outcomes, but
the number of trials available is limited by computational resources
or experimental effort. Changes in contact and discrete controller
decisions can further complicate the relationship between parameters
and performance. Effective tuning must therefore extract useful
information from a limited number of evaluations.

Uncertainty adds a further requirement: favorable performance on the
conditions used during tuning does not by itself establish reliable
behavior under different conditions. Performance criteria must account
for both efficient motion and successful task completion, and observed
improvements must be assessed across relevant operating conditions.
This creates two connected questions: how efficiently can a search
procedure identify promising parameter settings, and how consistently
does the resulting controller perform beyond the trials that guided
their selection?

\subsection{Related work}

To position the proposed method, we review population and subspace  optimization alongside contact-rich controller tuning and identify the combination studied here. We first discuss black-box/derivative-free optimization (BBO/DFO) methods,
with particular emphasis on population-based search and trust-region
strategies. We then examine subspace optimization techniques for constructing
and solving reduced-dimensional optimization models. Next, we review
learning-based methods and compliance control for contact-rich robotic
assembly, followed by derivative-free approaches to optimizing controller
parameters. Together, these research directions provide the theoretical
and methodological context for combining particle swarm optimization
with a reduced-subspace trust-region safeguard to tune an interpretable
robotic insertion controller.

\subsubsection{Derivative-free optimization}
\label{subsec:related-dfo}

BBO addresses problems for which candidate solutions
can be evaluated, but an explicit, computationally exploitable expression for
the objective or reliable derivative information is unavailable. Such problems
arise when evaluations involve numerical simulation, physical experiments, or
closed-loop control. DFO develops search
strategies that do not require the user to supply derivatives, encompassing
direct-search, model-based, stochastic, and population-based techniques.
General treatments are given by Audet and Hare~\cite{AudH}, Conn et
al.~\cite{conn2009}, and Larson et al.~\cite{LarMW}. Comparative studies by
Mor\'e and Wild~\cite{morewild2009} and Rios and
Sahinidis~\cite{RioS}, together with the methods described in
\cite{kimiaei2021vrbbo,VRDFON}, illustrate the diversity of these methods and the importance of evaluation-budget-aware benchmarking.

A further difficulty in BBO is that observed objective values may be
inaccurate or uncertain.
Numerical discretization and approximate models introduce deterministic
errors, whereas random simulation inputs, measurement errors, and stochastic
system dynamics can introduce variability. Consequently, an observed decrease
in the objective need not imply a comparable improvement in its underlying
expected value, especially when evaluations are expensive and few trials are
available. For the controller-tuning experiment considered here, a fixed
seed-specific batch of scenarios is used to compare candidates consistently
within each optimization run; screening and validation on independent
scenarios are then used to examine generalization beyond that batch.

Particle swarm optimization (\texttt{PSO}) is a population-based
DFO metaheuristic introduced by Kennedy and Eberhart~\cite{pso}. Each particle adjusts its velocity and position using
its own best previously observed solution and the best solution identified by
the swarm. This mechanism combines individual search experience with
population-level information and does not require a local model of the
objective. Its comparatively simple updates make it attractive for bounded
BBO problems, but the swarm may stagnate or devote a substantial
fraction of its evaluation budget to population updates that fail to improve
the incumbent. Moreover, the standard \texttt{PSO} update does not explicitly
construct a local approximation for systematic refinement of the current best
solution. These considerations motivate the addition of a conditional local
search mechanism, without replacing the underlying particle updates.

Vaz and Vicente~\cite{vazvicente2007pswarm} introduced
a particle swarm pattern search (\texttt{PSwarm}), which incorporates particle swarm exploration
into a pattern search framework for bound-constrained DFO.
When the swarm search fails to improve the incumbent,
the method polls coordinate directions around the current
best point. An unsuccessful poll contracts the mesh size,
allowing progressively finer local exploration.
Their pattern search analysis establishes stationarity results
under appropriate smoothness and mesh assumptions.
This provides an early example of combining population
exploration with conditional local refinement.

Trust-region (\texttt{TR}) methods instead seek improvement by approximately minimizing
a local model within a neighborhood of the current iterate. In DFO,
\texttt{TR} models can be built from sampled function values rather than
analytic gradients; trial steps are typically assessed by comparing
actual and predicted reductions~\cite{conn2009,LarMW}. Quadratic \texttt{TR}
subproblems can be characterized through optimality conditions and solved using
spectral or secular-equation techniques~\cite{moresorensen1983}. Reduced
subspaces can make these models less costly to store and manipulate in
higher dimensions~\cite{cartisroberts2023}. A model-based trial is not
automatically beneficial: the model may be inaccurate, its trial may violate
the bounds, and constructing it consumes evaluations. Therefore, a \texttt{TR}
safeguard within a population method should be activated selectively, and its
model probes and trial evaluations must share the overall evaluation budget. Several derivative-free \texttt{TR} solvers have been developed, including
\texttt{UOBYQA}~\cite{powell2002uobyqa},
\texttt{NEWUOA}~\cite{powell2006newuoa},
\texttt{BOBYQA}~\cite{powell2009bobyqa},
\texttt{BCDFO}~\cite{gratton2011bcdfo}, and, for nonsmooth objectives,
\texttt{DFO-TRNS}~\cite{liuzzi2019trustregion}.

Hybrid designs already demonstrate how DFO globalization
mechanisms can augment evolutionary search. The matrix adaptation derivative-free optimization (\texttt{MADFO}) method~\cite{kimiaei2024madfo} combines
matrix adaptation with a randomized nonmonotone line search for noisy
objectives. Its successor the matrix adaptation trust-region strategy (\texttt{MATRS})~\cite{kimiaei2025matrs} integrates mutation,
selection, recombination, trust-region steps, and a dedicated mixed-integer
search phase for noisy, bound-constrained mixed-integer problems. These methods
establish that a population-oriented search can be combined with additional
local-search or trust-region components, but their matrix-adaptation updates
and intended problem classes differ from those of the present work.

Together, these approaches highlight a recurring trade-off in
evaluation-limited DFO. Broad exploration can identify promising
regions, while local models use sampled information to guide more
focused improvement. The potential benefit of refinement must,
however, be weighed against the evaluations, storage, and numerical
work required to construct those models. This makes the choice of
model dimension and search directions important. Restricting a model
to a lower-dimensional space can reduce its computational cost, but
also limits the directions it represents. Subspace methods therefore
raise a further question: how can a compact representation retain
the information most useful for continued improvement?

\subsubsection{Subspace optimization methods}

Subspace optimization is a widely used strategy for addressing
high-dimensional problems by restricting search steps or local models to
lower-dimensional spaces. Its main potential benefit is to reduce the cost
of constructing and solving the relevant subproblems; a reduction in the
number of objective evaluations does not automatically follow.

In large-scale gradient-based optimization, limited-memory and active-set
methods exploit information on lower-dimensional spaces of search and
free-variable directions. For example, the limited-memory method for bound-constrained optimization (\texttt{LMBOPT}) algorithm of Kimiaei et al. \cite{kimiaei2022lmbopt} combines limited-memory
subspace models and search directions with an active-set strategy. Cartis et al. \cite{cartisfowkesshao2022} also developed a randomized
subspace framework for smooth nonconvex optimization, including
subspace Gauss--Newton methods for nonlinear least squares.

For nonlinear least-squares problems in DFO, Kimiaei and
Neumaier~\cite{kimiaei2022lmls} proposed the limited-memory method for nonlinear least squares (\texttt{LMLS}), a limited-memory trust-region method
that estimates residual Jacobian information using random coordinate
subspaces and incorporates nonmonotone acceptance and adaptive-radius
strategies. Cartis and Roberts developed the derivative-free block Gauss–Newton (\texttt{DFBGN}) framework, which constructs local interpolation models
within random subspaces to improve scalability on large nonlinear
least-squares problems~\cite{cartisroberts2023}.

For general DFO problems, Kimiaei, Neumaier, and
Faramarzi~\cite{kimiaei2023ssdfo} introduced the subspace method for unconstrained derivative-free optimization (\texttt{SSDFO}), which constructs quadratic models
in subspaces spanned by previous search directions and uses the resulting
models to obtain line-search directions.
More recently, Cartis and Roberts~\cite{cartisroberts2026quadratic} analyzed randomized subspace
quadratic models in DFO, including
complexity guarantees for approximate second-order stationarity.

Related work by Cartis, Massart, and Otemissov~\cite{cartis2023boundembeddings,cartis2023globalembeddings} investigates multiple
random embeddings for bound-constrained global optimization of functions
with low effective dimensionality and a broader random-embedding
framework for globally optimizing Lipschitz-continuous objectives.

These studies illustrate different ways of constructing reduced search
spaces: exploiting problem structure, sampling random directions, or
retaining information from previous iterations. Each choice emphasizes
different information and relies on particular assumptions.
The usefulness of a reduced model therefore depends not only on its
dimension, but also on whether the represented directions allow
meaningful improvement in the objective. For objectives defined by
robotic task performance, parameter changes act through feedback,
motion, and physical interaction. In contact-rich assembly, understanding
these effects requires attention to both the controller's behavior and
the uncertainty it encounters.

\subsubsection{Learning and compliance in assembly}
\label{subsec:learning-compliance-assembly}

Contact-rich insertion presents a control challenge because successful
assembly requires both accurate geometric alignment and appropriate
interaction with the environment. Uncertainty in the relative pose of
mating parts, friction, and clearance can turn an otherwise feasible
insertion motion into rim contact, jamming, or excessive loading.
Consequently, position tracking alone may be insufficient: the
controller must account for contact forces and determine when to
change its motion in response to stalled progress. Classical hybrid
position--force control separates position- and force-controlled task
directions~\cite{raibertcraig1981}, whereas impedance control addresses
the dynamic relationship between interaction forces and
motion~\cite{hogan1985impedance}. Both approaches provide important
conceptual foundations for contact-aware manipulation, although the
controller studied here does not implement the classical hybrid
position--force architecture or online impedance adaptation.

Learning-based methods seek to reduce the manual design and tuning
effort associated with contact-rich assembly. Inoue et
al.~\cite{inoue2017assembly} used reinforcement learning (RL) with
recurrent deep neural networks for high-precision peg-in-hole assembly and investigated robustness
to errors in position and orientation, including experiments on a
physical robot. Beltran-Hernandez et al.~\cite{compliance} studied
variable-compliance control for peg-in-hole insertion under uncertain
hole position, combining off-policy RL with domain
randomization and simulation-to-real transfer techniques. These
studies illustrate how learning can adjust motion and contact response
under geometric uncertainty rather than relying exclusively on a
fixed nominal trajectory.

A complementary approach preserves a conventional controller and
learns corrections to its output. Johannink et
al.~\cite{johannink2019residual} investigated residual RL,
superimposing a learned action on a hand-engineered feedback
controller and demonstrating the approach on a physical block-assembly
task. Zhang et al.~\cite{zhang2023} developed an offline--online
framework in which RL determines motion and
compliance-control parameters in randomized simulation, followed by
online residual adaptation of admittance parameters using force
measurements. Such approaches exploit existing control structure while
allowing the response to change under uncertain interaction
conditions; they differ from the offline fixed-parameter optimization
considered in this paper.

Simulation-based RL has also been developed for
multi-stage robotic assembly. \texttt{IndustReal}~\cite{industreal}
combines simulation-aware policy updates, signed-distance-field
rewards, and sampling-based curricula for contact-rich manipulation
and assembly; it also introduces a policy-level action integrator to
reduce execution errors during deployment. Its reported experiments
include transfer of learned policies to a physical robotic assembly
system. Force-guided exploration for robust contact-rich manipulation under uncertainty (\texttt{FORGE}) ~\cite{forge} addresses pose uncertainty through
force-aware policy learning, dynamics randomization, and conditioning
on an allowable force threshold. It additionally predicts task
success to support termination and force-threshold adjustment and
demonstrates physical contact-rich assembly tasks. Together, these
studies highlight the importance of evaluating contact-sensitive
policies beyond their training simulations.

These studies illustrate that learning can act at different levels
of contact-rich control, including task actions, corrections to an
existing feedback law, and motion or compliance parameters.
The choice of what to learn influences the search space and the
information needed to assess improvement. When a useful controller
structure is already available, selecting its numerical settings
still requires balancing progress, completion, and contact loading.
Observed task outcomes provide criteria for comparing candidate
settings, connecting contact-aware assembly with policy optimization
based on closed-loop performance. This connection motivates examining
how controller parameters can be selected when analytical performance
derivatives are unavailable.

\subsubsection{Derivative-free policy optimization}
\label{subsec:related-policy-optimization}

Derivative-free policy optimization (DFPO) treats the selection of feedback-law, motion-primitive, or staged-controller parameters as a BBO problem, using measured closed-loop performance over complete episodes.
In this setting, the objective may depend on contact events, switching
rules, termination conditions, or stochastic scenario realizations, and
its derivatives with respect to policy parameters may be unavailable or
unreliable. This form of direct policy search offers an alternative to
learning an unrestricted state-to-action mapping: it optimizes a
predefined parameterization using observed rollout costs or rewards.
In contrast, theoretical policy-optimization results for structured
linear--quadratic control may exploit differentiable policies and more
specific modeling assumptions~\cite{guo2026policy}. Such assumptions
are not imposed on the insertion-controller objective considered here.

Evolutionary strategies are established tools for episodic policy
search. The covariance matrix adaptation evolution strategy
(\texttt{CMAES})~\cite{hansen} adapts a multivariate sampling distribution over
continuous parameters, whereas the matrix adaptation
evolution strategy (\texttt{MAES})~\cite{maes} updates a sampling transformation
directly. In a recent robotic application,
Jafari-Tabrizi et al.~\cite{jafari2024inspection} used
\texttt{CMAES} and policy-improvement techniques to optimize
parameterized robot scanning trajectories from image-quality feedback.
Sago et al.~\cite{sago2025path} investigated multiobjective
path-integral policy improvement for robotic motion. These studies
illustrate how rollout-level performance measures can guide policy
parameter updates without requiring analytical differentiation through
the robot and its environment.

Population heuristics provide another route to controller tuning.
Blanck-Kahan et al.~\cite{blanck2024parallel} employed
\texttt{DE} to tune controller gains for a parallel robot, while
Rodr\'iguez-Molina et al.~\cite{rodriguez2024metaheuristic}
combined polynomial-regression identification with online
\texttt{DE}-based tuning for simulated robotic manipulators.
Hasan et al.~\cite{hasan2024spider} optimized a fuzzy mobile-robot
controller using a social-spider search strategy and compared its
simulated navigation performance with \texttt{PSO}-based alternatives.
Huang et al.~\cite{huang2025neurofuzzy} used \texttt{PSO}
and beetle-antennae search to tune a neuro-fuzzy active-force
controller in robotic trajectory-tracking simulations.
These examples motivate population-based parameter optimization, but
they do not establish that a particular metaheuristic is uniformly
best across different robotic tasks or evaluation budgets.

Data-efficient DFPO has also been investigated within model-based
and RL frameworks. Black-box data-efficient robot policy search 
(\texttt{Black-DROPS}) combines DFPO with an uncertain
learned dynamics model~\cite{blackdrops}.
Lei et al.~\cite{lei2022derivativefree} integrated DFO-based
policy updates into an actor--critic framework for automatic
controller tuning, reporting numerical studies of cruise control and
trajectory tracking. More recently, He et al.~\cite{he2025simulation}
combined robot and simulator observations in a probabilistic
model for DFPO to reduce required physical interactions.
Metaheuristic optimization has also been combined with value
approximation for nonlinear optimal control~\cite{amirabadi2024}.
These approaches differ in whether simulation data support model
learning, policy improvement, or direct optimization of a fixed
controller parameter vector.

Surrogate-assisted methods place particular emphasis on expensive
evaluations. Chen et al.~\cite{chen2024impedance} studied
multiobjective Bayesian optimization of robot impedance parameters
using transient and steady-state force-tracking objectives, with
experimental validation. Nobar et al.~\cite{nobar2025digitaltwin}
used guided Bayesian optimization and a data-driven digital twin for
controller tuning. For high-dimensional model-predictive-control
parameter selection, Kudva et al.~\cite{kudva2024sparse} proposed
Bayesian optimization with sparse-subspace surrogate models.
Zag\'orowska et al.~\cite{zagorowska2025safe} examined more
efficient safe learning for controller tuning and validated their
method on a precision motion system. These studies highlight the
benefits of learned performance models and explicit validation,
although constructing a global performance surrogate or providing
an in-experiment safety guarantee is not part of the present algorithm.

Across these approaches, controller tuning involves more than choosing
a search algorithm. The policy parameterization determines which
behaviors can change, while the objective and evaluation conditions
determine what counts as improvement. Evaluation efficiency must
therefore be interpreted alongside the quality and scope of the
evidence obtained about a controller. Meaningful assessment requires
explicit performance criteria, transparent accounting of evaluation
effort, and a distinction between improvement on the trials used for
parameter selection and performance under additional operating
conditions. These considerations provide a basis for formulating
tuning problems and assessing their solutions.

\subsection{Contributions and relation to existing methods}
\label{subsec:contributions}

The preceding review highlights three connected concerns: allocating
limited evaluations between exploration and local refinement,
selecting the information retained in reduced models, and assessing
controller performance beyond the trials used for parameter selection.
We investigate these concerns through a conditional integration of
\texttt{PSO} and reduced-subspace \texttt{TR} refinement, together with its
application to a predefined insertion-controller family. The numerical
assessment examines both benchmark optimization performance and the
complete controller-tuning and deployment pipeline.

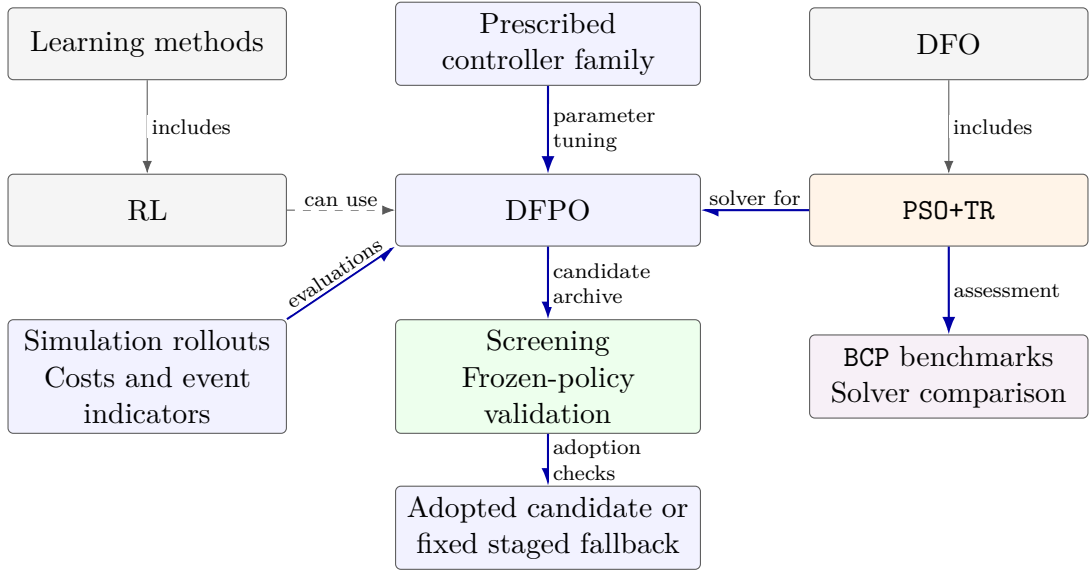
\begin{figure}[!htbp]
	\centering
	\begin{tikzpicture}[x=1cm,y=1cm,>=Latex,font=\small,
		base/.style={draw=black!65,rounded corners=2pt,
			text width=3.4cm,align=center,minimum height=.95cm,
			inner sep=4pt},
		context/.style={base,fill=gray!8},
		paper/.style={base,fill=blue!5,text width=3.75cm},
		solver/.style={base,fill=orange!9},
		check/.style={paper,fill=green!7},
		assessment/.style={base,fill=violet!6},
		flow/.style={->,draw=blue!65!black,thick},
		relation/.style={->,draw=black!65},
		optional/.style={relation,dashed},
		edge label/.style={font=\scriptsize,fill=white,inner sep=1.5pt}]
		
		\node[context] (learning) at (0,0) {Learning methods};
		\node[paper] (family) at (5.3,0) {Prescribed\\controller family};
		\node[context] (dfo) at (10.6,0)
		{DFO};
		
		\node[context] (rl) at (0,-2.2)
		{RL};
		\node[paper] (dfpo) at (5.3,-2.2)
		{DFPO};
		\node[solver] (psotr) at (10.6,-2.2)
		{\texttt{PSO+TR}};
		
		\node[paper,text width=3.4cm] (episodes) at (0,-4.4)
		{Simulation rollouts\\Costs and event indicators};
		\node[check] (screen) at (5.3,-4.4)
		{Screening\\Frozen-policy\\validation};
		\node[assessment] (bcp) at (10.6,-4.4)
		{\texttt{BCP} benchmarks\\Solver comparison};
		
		\node[paper] (deployed) at (5.3,-6.4)
		{Adopted candidate or\\fixed staged fallback};
		
		\draw[relation] (learning.south)
		--node[edge label,right] {includes} (rl.north);
		\draw[relation] (dfo.south)
		--node[edge label,right] {includes} (psotr.north);
		\draw[flow] (family.south)
		--node[edge label,right,text width=1.65cm] {parameter tuning}
		(dfpo.north);
		\draw[optional] (rl.east)
		--node[edge label,above] {can use} (dfpo.west);
		\draw[flow] (psotr.west)
		--node[edge label,above] {solver for} (dfpo.east);
		\draw[flow] (episodes.north east)
		--node[edge label,above,sloped] {evaluations} (dfpo.south west);
		\draw[flow] (dfpo.south)
		--node[edge label,right,text width=1.7cm] {candidate archive}
		(screen.north);
		\draw[flow] (screen.south)
		--node[edge label,right,text width=1.7cm] {adoption checks}
		(deployed.north);
		\draw[flow] (psotr.south)
		--node[edge label,right] {assessment} (bcp.north);
	\end{tikzpicture}
	\caption{High-level relationships between learning, RL, DFPO, and the \texttt{PSO+TR} solver in the present paper.}
	\label{fig:method-context}
\end{figure}

Figure~\ref{fig:method-context} distinguishes learning paradigms,
policy-optimization formulations, and numerical solvers.
RL is a learning paradigm, DFPO selects policy parameters from
episode evaluations without requiring objective derivatives, and
\texttt{PSO+TR} is the DFO solver used here.
RL can also use derivative-free policy search; the dashed connection
indicates this overlap rather than mutually exclusive categories.
The central route shows the prescribed controller family, parameter
tuning, and independent adoption checks used in this study.

The principal contributions are summarized as follows:
\begin{itemize}
	
	\item {\bf A conditional and successive \texttt{TR} safeguard for \texttt{PSO}.}
	The standard swarm-generation and personal/global-best updates are
	retained. A \texttt{TR} invocation is eligible only after a completed generation
	fails to improve the incumbent and sufficient evaluation allowances
	remain. Within an invocation, accepted steps may be followed by further
	reduced-space trials under explicit continuation rules.
	Initialization, population evaluations, model probes, and trial
	attempts share one charged evaluation budget.
	
	\item {\bf A persistent basis with incremental reduced-curvature updates.}
	An initially randomized orthonormal basis is retained across iterations
	and updated from successful population displacements, with random-column
	refreshes during prolonged non-improvement.
	Updates following accepted \texttt{TR} steps can change column signs but preserve
	the subspace span in exact arithmetic.
	A frozen swarm-geometry factor defines a positive-definite reduced
	curvature surrogate. Replacing one basis column requires updating only
	the corresponding factor column and reduced-matrix row and column,
	without constructing an ambient-dimensional curvature matrix.
	The surrogate represents search geometry and is not presented as an
	objective-Hessian estimate.
	
	\item {\bf A bounded DFPO formulation for contact-rich insertion.}
	We tune 14 parameters of a staged approach--search--alignment--insertion
	controller with history-triggered probing and retreat over
	$\Theta=[-3,3]^{14}$.
	The objective uses the larger expected penalized episode cost under
	nominal and development scenario distributions.
	The feedback law, phase logic, and recovery rules are prescribed;
	optimization selects their numerical parameters.
	The controller receives measured pose, velocity, and net contact wrench,
	without access to the true socket pose or contact parameters.
	
	\item {\bf Independent controller selection and evaluation.}
	Candidate screening uses scenarios excluded from optimization.
	The selected candidate is then frozen before independent validation,
	and adoption requires predefined screening and validation criteria.
	Otherwise, the fixed staged reference is retained.
	Final testing includes a withheld stress distribution that does not
	enter optimization or adoption decisions.
	This protocol separates the observations used to choose parameters
	from those used to assess the selected controller.
	
	\item {\bf Complementary numerical assessments.}
	We compare \texttt{PSO+TR} with four DFO solvers on 93 eligible
	\texttt{BCP} problems with dimensions 1--12 from the
	Branch-And-Reduce Optimization Navigator ({\tt BARON}) collection
	of Sahinidis~\cite{baron}, using common targets and matched evaluation
	budgets. Separately, eight independent optimization runs assess the
	controller-tuning pipeline on a reconstructed synthetic planar
	insertion model, with paired final-test comparisons against three
	fixed controllers.
	These experiments distinguish benchmark solver performance from
	application-specific controller outcomes.
	
\end{itemize}

\noindent{\bf Relation to optimization methods.}
Population methods such as \texttt{PSO}, \texttt{CMAES}, and
\texttt{MAES} provide established mechanisms for search without
objective derivatives~\cite{pso,hansen,maes}.
A close precedent is \texttt{PSwarm} ~\cite{vazvicente2007pswarm}, which combines swarm search
with coordinate polling when the search fails to improve the
incumbent. Hybrid methods such as \texttt{MADFO} and \texttt{MATRS} already
combine evolutionary search with additional local-search or
\texttt{TR} mechanisms~\cite{kimiaei2024madfo,kimiaei2025matrs}.
Likewise, \texttt{SSDFO} and \texttt{DFBGN} demonstrate the use of
reduced spaces for local optimization
models~\cite{cartisroberts2023,kimiaei2023ssdfo}.
The present design specifies how a completed-generation trigger,
successive \texttt{TR} attempts based on reduced quadratic models,
persistent basis information, and a frozen geometry factor
operate together within a \texttt{PSO} search.
Its methodological contribution concerns this integration and its
evaluation; none of the constituent mechanisms is claimed to be
new in isolation.

\noindent{\bf Relation to contact-rich policy learning and tuning.}
The use of motion and force feedback has classical foundations in
contact-aware control~\cite{raibertcraig1981,hogan1985impedance},
while direct optimization of parameterized robotic behavior is
established in controller and trajectory
tuning~\cite{jafari2024inspection,rodriguez2024metaheuristic}.
Learning-based assembly methods span policy learning, residual
correction, and compliance adaptation
\cite{industreal,compliance,forge,johannink2019residual}.
Here, the controller architecture is specified in advance and its
numerical parameters are selected offline.
The tuning procedure does not fit a dynamics model or a global
performance surrogate, as employed in model-learning and
surrogate-assisted approaches
\cite{blackdrops,he2025simulation,nobar2025digitaltwin};
its quadratic models serve local optimization.
The application-level contribution is the joint tuning of search,
insertion, and recovery parameters under uncertainty, coupled with
independent adoption checks and a specified fallback.
The robotic results assess this complete pipeline within the
synthetic simulator. A matched \texttt{PSO}-only robotic ablation
would be required to isolate the contribution of the \texttt{TR} safeguard.

\begin{figure}[!htbp]
	\centering
	\begin{tikzpicture}[x=1cm,y=1cm,>=Latex,font=\scriptsize,
		box/.style={draw=black!65,rounded corners=2pt,fill=blue!5,
			text width=5.9cm,align=center,
			minimum height=.62cm,inner sep=3pt},
		branch/.style={draw=black!65,rounded corners=2pt,fill=green!7,
			text width=2.85cm,align=center,
			minimum height=.67cm,inner sep=3pt},
		trial/.style={draw=black!65,rounded corners=2pt,fill=orange!9,
			text width=2.85cm,align=center,
			minimum height=.67cm,inner sep=3pt},
		edge label/.style={font=\scriptsize,fill=white,inner sep=1pt}]
		
		\node[font=\small\bfseries] at (2.1,3.82)
		{(a) Planar contact geometry};
		\fill[gray!23] (0.0,-0.15) rectangle (1.38,2.65);
		\fill[gray!23] (2.82,-0.15) rectangle (4.2,2.65);
		\draw[thick]
		(0,-0.15)--(0,2.65)--(1.38,2.65)--(1.38,-0.15);
		\draw[thick]
		(2.82,-0.15)--(2.82,2.65)--(4.2,2.65)--(4.2,-0.15);
		\node[rotate=90] at (.66,1.0) {Solid rim};
		\node[rotate=90] at (3.54,1.0) {Solid rim};
		\begin{scope}[shift={(2.1,2.12)},rotate=8]
			\filldraw[fill=cyan!17,draw=cyan!65!black,thick]
			(-.42,-.85) rectangle (.42,1.05);
		\end{scope}
		\draw[dashed,thin] (2.1,3.37)--(2.1,-.2);
		\draw[->,thick] (2.1,.85)--(2.1,.27);
		\node at (2.1,-.48) {Opening};
		
		\node[font=\small\bfseries] at (9.2,3.82)
		{(b) \texttt{PSO+TR} computational workflow};
		
		\node[box] (pso) at (9.2,3.05)
		{Complete \texttt{PSO} generation;
			update particles and bests};
		
		\node[box] (decision) at (9.2,2.08)
		{Did the generation improve the global best?};
		
		\node[branch] (yes) at (7.48,1.00)
		{Update subspace and continue \texttt{PSO}};
		
		\node[trial] (no) at (10.96,1.00)
		{Consider a \texttt{TR} invocation};
		
		\node[branch] (budget) at (10.96,-.35)
		{Do both evaluation budgets permit a \texttt{TR} attempt?};
		
		\node[trial] (trstep) at (10.96,-1.60)
		{Build reduced model; attempt \texttt{TR} step};
		
		\node[trial] (accepted) at (10.96,-2.95)
		{Accepted: update incumbent, basis,
			and reduced curvature};
		
		\node[branch] (failed) at (7.48,-2.95)
		{Earlier \texttt{TR} step\\accepted?};
		
		\node[branch] (retry) at (7.48,-4.50)
		{Minimum radius and both budgets permit
			a half-radius retry?};
		
		\node[trial] (halve) at (10.96,-4.50)
		{Retain center, basis, and curvature;
			halve working radius};
		
		\node[box] (exit) at (9.2,-6.00)
		{End \texttt{TR} with final accepted incumbent\\
			(or unchanged best if none).\\
			Resume \texttt{PSO} if the global budget permits.};
		
		\draw[->] (pso)--(decision);
		
		\draw[->] (decision.south)--++(0,-.15)
		-| node[pos=.75,edge label] {Yes} (yes.north);
		
		\draw[->] (decision.south)--++(0,-.15)
		-| node[pos=.75,edge label] {No} (no.north);
		
		\draw[->] (yes.west)--(5.65,1.00)|-(pso.west);
		\draw[->] (no.south)--(budget.north);
		
		\draw[->] (budget.south)
		--node[edge label,right] {Yes} (trstep.north);
		
		\draw[->] (budget.west)--(5.40,-.35)
		--(5.40,-6.00)--(exit.west);
		\node[edge label,above] at (7.15,-.35) {No};
		
		\draw[->] (trstep.south)--(accepted.north);
		
		\draw[->] (trstep.west)
		--node[edge label,above] {Unsuccessful} (7.48,-1.60)
		--(failed.north);
		
		\draw[->] (accepted.east)--++(.55,0)|-(budget.east);
		
		\draw[->] (failed.west)
		--node[edge label,above] {Yes} (5.40,-2.95);
		\fill (5.40,-2.95) circle (1.1pt);
		
		\draw[->] (failed.south)
		--node[edge label,right] {No} (retry.north);
		
		\draw[->] (retry.east)
		--node[edge label,above] {Yes} (halve.west);
		
		\draw[->] (retry.south)--(7.48,-5.20)
		-| (exit.north);
		\node[edge label,right] at (7.48,-5.30) {No};
		
		\draw[->] (halve.east)--(13.85,-4.50)
		|- (budget.north east);
		
		\node[edge label,rotate=90] at (14.3,-2.2) {
			\begin{tabular}{@{}l@{}}
				Before first success, if radius \\
				and budgets permit
			\end{tabular}
		};
		
	\end{tikzpicture}
	
	\caption{Illustrative planar peg--rim geometry and conditional \texttt{PSO+TR} workflow.}
	\label{fig:workflow}
\end{figure}
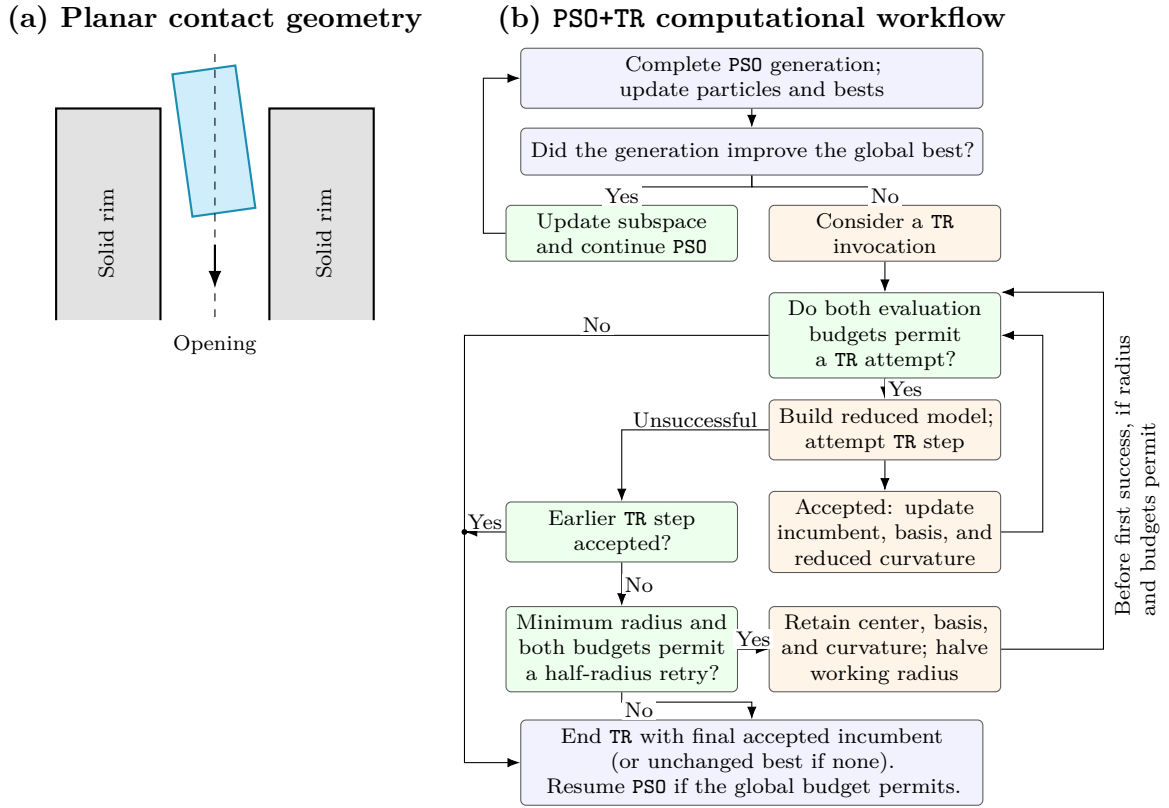

Figure~\ref{fig:workflow} combines the schematic planar peg--rim
geometry with the conditional search operations of \texttt{PSO+TR}.
The left panel is illustrative, not a measured trajectory.
After a complete non-improving \texttt{PSO} generation, a \texttt{TR}
invocation may begin if both the global and per-invocation evaluation
allowances permit. Before the first accepted \texttt{TR} step, an
unsuccessful attempt may be retried at half the working radius if
both allowances and the minimum radius permit. Accepted steps update
the incumbent and reduced state and may be followed by further
\texttt{TR} attempts. After any accepted step, the first subsequent
unsuccessful attempt ends the invocation. The invocation also ends
when the evaluation allowances do not permit continuation.
Control returns to \texttt{PSO} if the global budget permits.
Independent screening, frozen-policy validation, and paired final
evaluation are described in Section~\ref{sec:robotic-selection}.

\clearpage

\section{Tuning problem}
\label{sec:tuning-problem}

To balance insertion cost and completion reliability under uncertainty,
we formulate a DFPO problem over 14 parameters of a structured insertion
controller. We minimize the larger of the expected penalized episode
costs under the nominal and development scenario distributions,
with fixed penalties for unsuccessful completion and benchmark-limit
violations. The controller's prescribed phase architecture, feedback
structure, and recovery rules remain fixed, while its numerical
parameters are tuned using closed-loop simulation outcomes.
This yields a finite-dimensional BBO problem, whose sample-based
approximation we optimize using the proposed \texttt{PSO+TR} algorithm.

The feasible normalized parameter set for the reported experiment is
$\Theta=[-3,3]^{14}\subset\R^{14}$.
Let $\theta\in\Theta$ denote the normalized controller parameter vector,
and let $\pi_\theta$ represent the corresponding insertion policy.
The components of $\theta$ determine the motion amplitudes,
controller gains, timing parameters, and recovery thresholds.

Under the physical parameter map
$p_j=c_j+r_j\tanh(\theta_j)$, $j=1,\ldots,14$,
the centers $c_j$ and positive radii $r_j$ are specified in
Tables~\ref{tab:policy} and~\ref{tab:hybridparams}, below.
Consequently, each mapped parameter satisfies
\[
p_j\in[c_j-r_j\tanh(3),\,c_j+r_j\tanh(3)],
\qquad j=1,\ldots,14,
\]
where $\tanh(3)\simeq0.995054754$.
The endpoints of these intervals are attained at
$\theta_j=\pm3$. The values $c_j\pm r_j$ are asymptotic
limits of the mapping and are not attained.

For a scenario $\omega$, let $C(\pi_\theta,\omega)$ denote the
episode cost, and define
$S(\pi_\theta,\omega)\in\{0,1\}$ and
$R(\pi_\theta,\omega)\in\{0,1\}$ as the completion and
benchmark-limit violation indicators, respectively.
Completion and violation are recorded as separate outcomes:
a successful episode may still involve a benchmark-limit violation.

Let $P_0$ and $P_1$ denote the nominal and development scenario
distributions, respectively. To account for the different
operating conditions, we formulate the bounded tuning problem as
\begin{equation}\label{eq:problem}
	\min_{\theta\in\Theta}
	\max_{j\in\{0,1\}}
	\E_{\omega\sim P_j}
	\left[
	C(\pi_\theta,\omega)
	+10\bigl(1-S(\pi_\theta,\omega)\bigr)
	+20R(\pi_\theta,\omega)
	\right].
\end{equation}

The objective combines episode cost, completion reliability,
and benchmark-limit violations. Taking the maximum of the
expected penalized costs over the two distributions discourages
solutions that perform well under nominal conditions but poorly
under the development distribution, or vice versa. Thus, the
formulation seeks distribution-level robustness with respect
to these two prescribed scenario families. It does not represent
worst-case optimization over every possible contact configuration,
nor does it guarantee satisfactory performance on distributions
that are excluded from the tuning objective.

The coefficients $10$ and $20$ are fixed penalties associated
with unsuccessful completion and benchmark-limit violations,
respectively. Their purpose is to discourage solutions that
achieve lower episode costs through unreliable insertion or
excessive contact loading. These coefficients represent
benchmark-specific preferences rather than physical actuator
limits or formally calibrated safety requirements.

The penalties in~\eqref{eq:problem} are additional to those already
included in the episode cost $C$, defined in~\eqref{eq:cost}.
All penalty weights remain fixed throughout the reported experiments;
no independent penalty-weight sensitivity study is claimed.

Evaluating the objective requires executing the parameterized
controller under uncertain contact conditions. Its dependence
on nonlinear contact dynamics, discrete phase transitions,
recovery decisions, and event-dependent termination makes
analytical differentiation difficult. We therefore treat
the objective as a BBO oracle and solve its sample-based
approximation using \texttt{PSO+TR}. Within each independent
optimization run, the same seed-specific scenario batch is
used for candidate evaluations, including population updates
and trust-region model construction. This permits consistent
comparisons of observed objective values, while the resulting
sample-average objective remains an approximation of the
distribution-level formulation.

Finally, the expectation-based objective and its soft penalties
do not impose hard probabilistic safety constraints. Accordingly,
optimization is followed by independent candidate screening
and frozen-policy validation, using scenarios excluded from
the optimization batch. These procedures assess completion,
violation, and paired-loss behavior before deciding whether
a tuned controller is adopted. Final evaluation on additional
nominal, development, and withheld stress scenarios provides
a separate assessment of the resulting controller within the
specified numerical benchmark.

\section{Planar insertion model and controller family}\label{sec:PlanarInsertion}

The DFPO formulation in Section~\ref{sec:tuning-problem} requires a
closed-loop simulation model for evaluating candidate controller parameters.
We therefore consider a synthetic planar insertion benchmark that couples
a four-corner compliant contact model with a structured, interpretable
controller. This section specifies the dynamics and scenario uncertainty,
the information available to the controller, its 14-parameter family,
and the episode costs and outcome indicators used by the BBO oracle.
The model provides the numerical setting for the subsequent
\texttt{PSO+TR} experiments; it is not a model of a complete physical robot.

\subsection{Dynamics, information, and uncertainty}

The planar configuration is $q=(x,z,\varphi)$, where $x$ and $z$
denote the peg-center position and $\varphi$ its orientation;
the positive $z$ direction points downward. The dynamics are
\begin{equation}
	\operatorname{diag}(m,m,I)\ddot q=u+w_c(q,\dot q;\omega).
\end{equation}
Here $u$ is the applied control wrench and
$w_c(q,\dot q;\omega)$ is the contact wrench for scenario $\omega$.
The mass is $m=0.12$ kg, inertia $I=2\times10^{-5}$ kg m$^2$,
peg half-width $a=5$ mm, and half-length $b=12$ mm.
The initial center is $(0,-17\text{ mm})$ with zero orientation
and velocity. These are synthetic benchmark parameters.
The model represents local planar motion; articulated joints,
gravity, perception, and actuator dynamics are not represented.
At each physics step, semi-implicit Euler integration updates velocity
using the current control and contact wrenches and then updates position
using the new velocity. The contact wrench is subsequently recomputed.
The control command is held constant between 10~ms controller updates.

An episode scenario $\omega=(d_x,\alpha,\mu,g,k_c)$ specifies socket
lateral position, orientation, friction, per-side clearance, and
contact stiffness. Its components are independently sampled from
uniform distributions with ranges in Table~\ref{tab:uncertainty}.
The nominal and development distributions are used to construct
the tuning objective in~\eqref{eq:problem}; the stress distribution
is withheld from tuning.
The controller receives measured pose, velocity, net contact wrench,
and elapsed time, but not the true socket pose, friction coefficient,
clearance, or stiffness. These scenario parameters are used internally
by the simulator and for outcome assessment. The DFO algorithm receives
only the resulting episode evaluations, not the hidden physical parameters.
This reconstructed \texttt{MATLAB} implementation supplies the BBO oracle
for the \texttt{PSO+TR} controller-tuning experiments reported in
Section~\ref{sec:psotr-robot-results}.

\begin{table}[ht]\centering\small
	\caption{Independent uniform scenario distributions.
		Nominal and development scenarios are used for tuning;
		stress scenarios are withheld for evaluation.}
	\label{tab:uncertainty}
	\begin{tabular}{llll}
		\toprule
		Parameter & Nominal & Development & Stress\\
		\midrule
		$d_x$ (mm) & $[-1.5,1.5]$ & $[-2.25,2.25]$ & $[-3,3]$\\
		$\alpha$ (rad) & $[-.04,.04]$ & $[-.06,.06]$ & $[-.08,.08]$\\
		$\mu$ & $[.12,.35]$ & $[.35,.55]$ & $[.55,.80]$\\
		$g$ (mm) & $[.35,.70]$ & $[.20,.35]$ & $[.10,.20]$\\
		$k_c$ (N/m) & $[1200,2000]$ & $[1200,2000]$ & $[1800,3000]$\\
		\bottomrule
	\end{tabular}
\end{table}

\subsection{Nearest-face contact}
The contact model determines the wrench $w_c$ appearing in the planar
dynamics. Let
\[
Q(\psi)=\begin{bmatrix}
	\cos\psi&-\sin\psi\\
	\sin\psi&\cos\psi
\end{bmatrix},
\qquad
J=\begin{bmatrix}0&-1\\1&0\end{bmatrix}.
\]
Index the four body-frame corners by $i=1,\ldots,4$, with
\[
\eta_i\in\{(-a,-b)^\top,(a,-b)^\top,(-a,b)^\top,(a,b)^\top\}.
\]
Their coordinates in the true socket frame are
\[
\begin{pmatrix}\xi_i\\\zeta_i\end{pmatrix}
=Q(\alpha)^\top\left[
\begin{pmatrix}x-d_x\\z\end{pmatrix}+Q(\varphi)\eta_i\right].
\]
For the stationary socket, the socket-frame moment arm measured
from the peg center and the corner velocity are
\[
r_i=Q(\alpha)^\top Q(\varphi)\eta_i,\qquad
v_i=Q(\alpha)^\top\left[
\begin{pmatrix}\dot x\\\dot z\end{pmatrix}
+\dot\varphi JQ(\varphi)\eta_i\right].
\]
All corner normals, velocities, and forces below are expressed
in the socket frame.

Let $p_i=|\xi_i|-a-g$. A corner penetrates a solid rim quadrant when
$p_i>0$ and $\zeta_i>0$. The scalar penetration depth is
\begin{equation}
	d_i=\max\{\min(p_i,\zeta_i),0\}.
\end{equation}
For active contact, the contact normal used by the force law is
$n_i=(-\operatorname{sgn}\xi_i,0)^\top$ if $p_i<\zeta_i$,
and $(0,-1)^\top$ otherwise.

Writing $n_i=(n_{x,i},n_{z,i})^\top$, define the tangent
$t_i=(-n_{z,i},n_{x,i})^\top$ and velocity components
$v_{n,i}=n_i^\top v_i$, $v_{t,i}=t_i^\top v_i$.
Throughout the paper, $[s]_+=(s)_+=\max\{s,0\}$.
For active contact, the nonnegative normal-force magnitude
$N_i$ and corner force $f_i$ are
\begin{equation}
	N_i=[k_cd_i-c_nv_{n,i}]_+,\qquad
	f_i=N_in_i-\mu N_i\tanh(v_{t,i}/v_0)t_i,
\end{equation}
where $c_n=1.5$ N s/m and $v_0=0.001$ m/s.
Inactive corner forces are zero.

Writing $r_i=(r_{x,i},r_{z,i})^\top$ and
$f_i=(f_{x,i},f_{z,i})^\top$, the world-frame contact wrench is
\[
w_c=
\begin{pmatrix}
	Q(\alpha)\sum_{i=1}^4 f_i\\
	\tau_c
\end{pmatrix},
\qquad
\tau_c=\sum_{i=1}^4
\left(r_{x,i}f_{z,i}-r_{z,i}f_{x,i}\right).
\]
Here, $\tau_c$ is the planar contact torque about the peg center.
The monitored load $L=\sum_i\|f_i\|_2$ sums the magnitudes of all
corner forces and therefore generally differs from the norm of the
net contact force supplied to the controller. This distinction matters
because the load-based violation threshold uses $L$, whereas controller
transitions use components of the net contact wrench.

The nearest-face rule avoids penalizing a side-contacting corner according to
its entire depth below the top plane. The rule remains a four-point compliant
approximation, with discontinuous normal changes at the nearest-face switch.
It does not resolve complete edge/face contact or three-dimensional geometry.

\subsection{Controller structure and parameterization}
The controller family $\pi_\theta$ introduced in
Section~\ref{sec:tuning-problem} comprises approach, search,
alignment, and insertion phases. Its structure and switching rules
remain fixed while DFO selects the numerical parameters.
Approach enters search when measured axial contact force is below $-0.10$ N
or measured center depth exceeds $-12$ mm. Search enters alignment at depth
$-9.5$ mm, latching the measured lateral coordinate and angle. Alignment enters
insertion after an episode-fixed dwell parameter and depth exceeding $-8$ mm.
The phase transitions use observable signals, without exact socket pose.

During approach and search, the lateral reference is
$A\sin(2\pi(t-t_s)_+/P)$ and the angular reference is zero,
where $t_s$ denotes the reference-time origin; this origin is reset
at the beginning of search. The oscillatory search motion can
therefore also act during approach. Alignment holds its entry coordinates.
Insertion sets lateral and angular targets to current measured coordinates,
retaining damping while releasing their static position errors.
The axial targets for the four phases are
$\min(-.017+v_at,-.010)$, $-.008$, $-.005$, and $.016+\Delta z$ m.

Each axis uses proportional--derivative (PD) wrench control,
\begin{equation}
	u_j=\operatorname{clip}\left(K_j(q_{d,j}-q_j)
	-2\sqrt{m_jK_j}\dot q_j,-\bar u_j,\bar u_j\right).
\end{equation}
Here, $j\in\{1,2,3\}$ indexes $x,z,\varphi$, respectively;
$q_{d,j}$ is the phase-dependent desired coordinate, and
$K_j$ is its proportional gain, as specified in Table~\ref{tab:policy}.
The effective inertias satisfy $(m_1,m_2,m_3)=(m,m,I)$.
The clipping operator is
$\operatorname{clip}(a,l,h)=\min\{h,\max\{l,a\}\}$.
The wrench caps are
$(\bar u_1,\bar u_2,\bar u_3)=(0.7,0.6,0.012)$
in N, N, and N m, respectively.
Gains and motion parameters remain fixed within an episode, except for the
specified phase-dependent switch. The isolated-axis damping formula does not
establish stability of the switched contact loop.

The first ten components of $\theta$ determine bounded physical
parameters through $p_j=c_j+r_j\tanh(\theta_j)$, $j=1,\ldots,10$.
The row order of Table~\ref{tab:policy} defines the coordinate order
of $\theta_j$ and the corresponding physical parameter $p_j$
for these ten components.
The four rows of Table~\ref{tab:hybridparams} correspond, in order,
to $j=11,\ldots,14$ and use the same mapping.
In both tables, each row specifies the center $c_j$ and radius $r_j>0$
in the listed units; values in mm and mm/s are converted to
m and m/s, respectively, for simulation.
Setting the first ten components to zero gives the nominal parameter
values $p_j=c_j$; disabling the additional probing/recovery guards
then recovers the fixed staged reference controller.

\begin{table}[ht]\centering\small
	\caption{Staged-policy parameter map.}\label{tab:policy}
	\begin{tabular}{lrrl}\toprule
		Parameter & Center $c_j$ & Radius $r_j$ & Unit\\
		\midrule
		Search amplitude $A$ & 3.5 & 2.0 & mm\\
		Search period $P$ & 4 & 2 & s\\
		Approach speed $v_a$ & 10 & 4 & mm/s\\
		Lateral gain before alignment & 70 & 30 & N/m\\
		Lateral gain after alignment & 12 & 6 & N/m\\
		Axial gain & 65 & 25 & N/m\\
		Angular gain before alignment & .025 & .0125 & N m/rad\\
		Angular gain after alignment & .006 & .003 & N m/rad\\
		Alignment dwell & .20 & .10 & s\\
		Final axial offset $\Delta z$ & 1.0 & .5 & mm\\
		\bottomrule
	\end{tabular}
\end{table}

\subsection{Objective and outcome definitions}
This subsection makes the cost and indicators in
\eqref{eq:problem} explicit. Let $T_e$ be the first recorded
completion or protective-termination time, capped at the horizon
of $8$ s. If neither event occurs, the episode ends at $T_e=8$ s.
The finite-horizon, undiscounted cost is
\begin{equation}\label{eq:cost}
	C=\int_0^{T_e}\left(.3+\frac{(.016-z)_+}{.026}
	+.02\left\|D_u^{-1}u\right\|^2\right)dt
	+12I_{\rm abort}+6I_{\rm unsuccessful}
	+2\frac{(.016-z(T_e))_+}{.026}.
\end{equation}
The indicator $I_{\rm abort}$ equals one for protective termination
and zero otherwise, and $I_{\rm unsuccessful}=1-S$.
A timeout without completion is unsuccessful even if no protective
abort occurs. Since $I_{\rm unsuccessful}=1-S$, the additional failure penalty
in~\eqref{eq:problem} combines with the penalty already included
in $C$ to give a total explicit unsuccessful-completion penalty
of $16(1-S)$. The protective-abort penalty $12I_{\rm abort}$
and violation penalty $20R$ are added separately; the violation
penalty applies regardless of the completion outcome. The matrix $D_u=\operatorname{diag}(3,3,.03)$ supplies fixed
effort reference scales in N, N, and N m, distinct from the command caps.
The composite cost is a benchmark objective, not a direct measurement
of physical energy. The progress terms in $C$ use world-frame depth,
whereas the completion indicator is evaluated in the true socket frame.
The control period is 10 ms.

For a stationary socket centered laterally at $d_x$ and oriented
by $\alpha$, define the socket-relative pose by
\[
\begin{pmatrix}x_{\rm rel}\\z_{\rm rel}\end{pmatrix}
=Q(\alpha)^\top\begin{pmatrix}x-d_x\\z\end{pmatrix},
\qquad
\varphi_{\rm rel}=\varphi-\alpha.
\]
The corresponding axial speed is
\[
\dot z_{\rm rel}=-\sin\alpha\,\dot x+\cos\alpha\,\dot z.
\]

The completion indicator $S$ equals one if completion is recorded
and zero otherwise. Completion requires all the following conditions
to hold simultaneously for 0.10 s at the physics sampling rate:
socket-frame center depth $z_{\rm rel}\geq16$ mm,
axial speed $|\dot z_{\rm rel}|<0.02$ m/s, and
\begin{equation}
	|x_{\rm rel}|+a|\cos\varphi_{\rm rel}|+b|\sin\varphi_{\rm rel}|
	\leq a+g+20\times10^{-6}\text{ m}.
\end{equation}
The full peg must also be below the rim:
$z_{\rm rel}\geq a|\sin\varphi_{\rm rel}|+b|\cos\varphi_{\rm rel}|$.
The 20 $\mu$m allowance accommodates small compliant penetration.

The violation indicator $R$ records any $L>1.5$ N,
$|\tau_c|>0.12$ N m, or protective workspace termination.
Protection terminates at $L>12$ N, $|\tau_c|>0.18$ N m,
$|x|>25$ mm, $|z|>50$ mm, or $|\varphi|>0.7$ rad.
Success and violation are separate labels.
State and accumulated cost freeze after termination.
Physics-rate checks cannot exclude an unobserved inter-step contact peak.
These limits specify the numerical benchmark and are not hardware ratings.

\subsection{History-triggered probing and retreat}
To recover from stalled insertion, the controller supplements its
staged motion with three insertion behaviors: release (mode 0),
local probing (mode 1), and retreat (mode 2). The mode numbers
identify controller states, not privileged physical contact states;
in particular, they are distinct from the peg mass $m$.
At each 10 ms control update, an exponentially weighted axial
velocity is updated by
\begin{equation}
	\bar v_{z,k}=\beta\bar v_{z,k-1}+(1-\beta)v_{z,k},
	\qquad \beta=0.935506985.
\end{equation}
Here, $k=0,1,\ldots$ indexes control updates at $t_k=k\Delta t_c$,
with $\Delta t_c=0.01$ s. The quantity $v_{z,k}=\dot z(t_k)$
is the measured world-frame axial velocity, and
$\bar v_{z,-1}=0$ initializes the filter at the start of each episode.
The decay time is $-\Delta t_c/\log\beta\simeq0.15$ s.

During insertion, a stall age accumulates when
$\bar v_z<v_{\rm stall}$, measured axial wrench is below $-0.10$ N,
and world center depth is less than 17 mm; otherwise it resets.
Release enters probing when $\text{stall age}\geq t_{\rm probe}$. Upon entering probing, the controller records the measured lateral position $x_{\rm entry}$ and entry time $t_{\rm entry}$.
The probe target is
\begin{equation}
	x_d=x_{\rm entry}+A_{\rm probe}\sin\!\left(
	2\pi(t-t_{\rm entry})/P_{\rm probe}\right),
	\qquad \varphi_d=0,
\end{equation}
where $P_{\rm probe}=0.6\ \mathrm{s}$ is the fixed probe period.
The existing post-alignment gains are retained.
During probing, release resumes when $\bar v_z\geq2v_{\rm stall}$,
and stall age resets.
A persistent stall triggers retreat when
$\text{stall age}\geq t_{\rm retreat}$, provided fewer than three
retreats have been completed. The target depth is latched
1.5 mm above the current measured depth.
The retreat-age accumulator is reset to zero on entry and
incremented by $\Delta t_c=0.01$ s at every retreat control
update, including the entry update. Search restarts when
$\text{retreat age}\geq0.25\ \mathrm{s}$. The search-to-alignment depth guard is re-evaluated at the next
10 ms control update. If the measured world-frame center depth
satisfies $z>-9.5\ \mathrm{mm}$ at that update, the controller
enters alignment. At most three retreats are allowed; later persistent stalls
remain in probing. Command caps are unchanged. The guards react to observed wrench and motion and do not receive
friction, clearance, or true socket pose.

Four additional components of $\theta$ parameterize the stall,
probe, and retreat responses using the same bounded mapping
$p_j=c_j+r_j\tanh(\theta_j)$
(Table~\ref{tab:hybridparams}), giving 14 tuning parameters in all.

\begin{table}[ht]\centering\small
	\caption{Additional hybrid-policy parameters. Total search dimension is 14.}
	\label{tab:hybridparams}
	\begin{tabular}{lrrl}\toprule
		Parameter & Center & Radius & Unit\\
		\midrule
		$v_{\rm stall}$ & 1.0 & .9 & mm/s\\
		$t_{\rm probe}$ & .30 & .15 & s\\
		$A_{\rm probe}$ & 1.5 & 1.0 & mm\\
		$t_{\rm retreat}$ & .90 & .40 & s\\
		\bottomrule
	\end{tabular}
\end{table}

For comparison, the fixed staged baseline disables probing and retreat
and uses zero values for the first ten normalized parameters.
The fixed-hybrid comparator enables these behaviors with all normalized
parameters set to zero. Both are predefined reference controllers,
not optimal recovery policies. Hence, the DFPO task is to choose
the parameters of the structured hybrid controller, not its phase
architecture or the contact geometry. The peg retains the noncircular,
rectangular geometry specified above.

\section{\texttt{PSO+TR} with a reduced-subspace \texttt{TR} safeguard}
\label{sec:pso-tr}

To supplement population exploration with local refinement,
\texttt{PSO+TR} combines \texttt{PSO} with a conditional
reduced-subspace \texttt{TR} procedure. A \texttt{TR} invocation
is eligible only after a complete \texttt{PSO} generation fails
to improve the incumbent and the remaining evaluation allowances
permit an attempt.

The working \texttt{TR} radius bounds the trial displacement and
the region in which the local model is used. Before the first
accepted step, an unsuccessful attempt may be retried at half
the working radius if the minimum radius and evaluation allowances
permit. Trial acceptance requires a strict objective decrease and
an actual-to-predicted reduction ratio meeting the prescribed
threshold. Accepted steps may be followed by further local
refinement; the first subsequent unsuccessful attempt ends
the invocation, as does exhaustion of either evaluation allowance.

The general formulation is the bounded DFO problem
\begin{equation}\label{eq:pso-problem}
	\min_{x\in[\ell,u]} f(x),\qquad x\in\R^n,
\end{equation}
where $f$ is evaluated through an oracle without using its
derivatives, and $\ell,u\in\R^n$ are componentwise bounds.

In the robotic application, $x$ represents the normalized controller
parameter vector $\theta$ from~\eqref{eq:problem}, with $n=14$ and
feasible set $\Theta=[-3,3]^{14}$. Accordingly, the bounds are
$\ell=-3\mathbf1_{14}$ and $u=3\mathbf1_{14}$,
where $\mathbf1_{14}\in\R^{14}$ is the all-ones column vector.
The same solver is evaluated on eligible \texttt{BCP} problems
through the five-method comparison in
Section~\ref{sec:bcp-profiles}, separately from its application
to robotic controller tuning.

A complete \texttt{PSO} generation precedes each possible \texttt{TR}
intervention. If the population improves its global best, the
safeguard is skipped. Otherwise, an eligible intervention can solve
multiple reduced \texttt{TR} subproblems. After an accepted step, the next
model is constructed around the newly accepted point and the process
continues while its evaluation allowance permits. A failure following
at least one success ends the intervention; unsuccessful trials before
the first success may instead be retried at smaller working radii.
This successive-step mechanism is distinct from repeatedly invoking
\texttt{TR} after every particle update.

For DFPO, the distribution-level objective in~\eqref{eq:problem}
is replaced during each optimization run by a sample-based BBO oracle.
Feasible population, \texttt{TR} probe, and \texttt{TR} trial points
are evaluated on the same fixed, seed-specific scenario batch. Thus, their objective
values can be compared consistently within a run, although observed
sample-average improvements need not generalize to new scenarios.
Independent screening and validation are described in
Section~\ref{sec:robotic-selection}.

\subsection{\texttt{PSO} population dynamics and conditional intervention}

Let $N$ be the swarm size. At generation $g=0,1,\ldots$,
particle $i\in\{1,\ldots,N\}$ has position $x_{i,g}\in\R^n$,
velocity $v_{i,g}\in\R^n$, and personal best $p_{i,g}\in\R^n$;
the current swarm best is $g_g\in\R^n$.
The symbol $\odot$ denotes componentwise multiplication, and
$\mathbf1_n\in\R^n$ is the all-ones vector.

For a finite box $[\ell,u]$, draw mutually independent vectors
$r_i^x,r_i^v\sim\mathcal U([0,1]^n)$ and initialize
\[
x_{i,0}=\ell+(u-\ell)\odot r_i^x,\qquad
v_{i,0}=0.1(u-\ell)\odot(2r_i^v-\mathbf1_n),
\qquad i=1,\ldots,N.
\]
Replace $x_{1,0}$ by the prescribed initial candidate
$x^0\in[\ell,u]$, evaluate all $N$ initial particles, and set
\[
p_{i,0}=x_{i,0},\qquad i=1,\ldots,N,
\qquad
g_0\in\argminx_{x\in\{p_{i,0}:1\le i\le N\}} f(x).
\]
All initialization evaluations are charged to the global
evaluation budget. The robotic runs use $N=10$ and
$x^0=\mathbf0_{14}$, the zero normalized parameter vector.


With independent vectors $r_{1,i,g},r_{2,i,g}$ uniformly distributed
on $[0,1]^n$, the \texttt{PSO} update~\cite{pso} is
\begin{align}
	v_{i,g+1}&=\omega v_{i,g}
	+c_1r_{1,i,g}\odot(p_{i,g}-x_{i,g})
	+c_2r_{2,i,g}\odot(g_g-x_{i,g}),
	\label{eq:pso-velocity}\\
	x_{i,g+1}&=\Pi_{[\ell,u]}(x_{i,g}+v_{i,g+1}),
	\label{eq:pso-position}
\end{align}
where $\Pi_{[\ell,u]}$ is the componentwise box projection.
The inertia weight is $\omega$, while $c_1$ and $c_2$ are the
cognitive and social acceleration coefficients, respectively.
The executed configuration uses $(\omega,c_1,c_2)=(0.7,1.4,1.4)$;
this scalar inertia weight is distinct from the scenario vector
$\omega$ used in the robotic model. When projection changes a
particle coordinate, the corresponding velocity component is set
to zero. These initialization, update, and boundary-handling rules
concern the proposed solver; the separate \texttt{PSO} comparator
retains its own implementation.

Each updated particle is evaluated before the personal and global
bests are recomputed:
\begin{equation}\label{eq:pso-pbest}
	\begin{aligned}
		p_{i,g+1}&=\begin{cases}
			x_{i,g+1},& f(x_{i,g+1})<f(p_{i,g}),\\
			p_{i,g},&\text{otherwise},
		\end{cases}\\
		g_{g+1}^{\rm PSO}
		&\in\argminx_{x\in\{p_{i,g+1}:1\le i\le N\}} f(x).
	\end{aligned}
\end{equation}
For robotic tuning, let $\Omega_j^{(s)}$ be the fixed scenario batch
for distribution $P_j$ and optimization seed $s$. The oracle is
\begin{equation}\label{eq:pso-score}
	\widehat f^{(s)}(\theta)
	=\max_{j\in\{0,1\}}
	\frac{1}{|\Omega_j^{(s)}|}
	\sum_{\xi\in\Omega_j^{(s)}}
	\left[C(\pi_\theta,\xi)
	+10\bigl(1-S(\pi_\theta,\xi)\bigr)
	+20R(\pi_\theta,\xi)\right].
\end{equation}
The robotic runs use $|\Omega_0^{(s)}|=|\Omega_1^{(s)}|=8$;
the training batches differ across seeds but remain unchanged
between evaluations within each run.

Let $f_g^{\rm old}=f(g_g)$ be the incumbent value immediately
before generation $g+1$. The nonnegative $\varepsilon_{\rm imp}$
is an absolute improvement tolerance in objective units; it is
zero in all reported runs. The reduced dimension used below is
$m=\min\{n,10\}$. A completed generation improves the incumbent when
\begin{equation}\label{eq:pso-trigger}
	f(g_{g+1}^{\rm PSO})
	<f_g^{\rm old}-\varepsilon_{\rm imp}.
\end{equation}
In that case, the best point is retained, the successful displacement
updates the persistent subspace, and no \texttt{TR} intervention occurs.
Otherwise, \texttt{TR} may be invoked only after every particle in that
generation has been evaluated and at least $m+1$ objective calls
remain in both the global budget and the intervention allowance.
There is at most one {\bf invocation} per eligible generation,
but that invocation may contain several reduced \texttt{TR} subproblems.

\subsection{Initialization and success-driven update of the search subspace}
\label{subsec:pso-subspace}

The reduced models are motivated by subspace DFO methods
\cite{cartisroberts2023,cartisroberts2026quadratic}.
Our implementation retains an orthonormal basis across iterations
and modifies it using successful displacement directions rather
than independently resampling an entire subspace at each call.
Set $m=\min\{n,10\}$ and, once per run, generate
\begin{equation}\label{eq:subspace-init}
	G_0=(\gamma_{ij})\in\R^{n\times m},\qquad
	\gamma_{ij}\overset{\mathrm{iid}}\sim\mathcal N(0,1),\qquad
	G_0=S_0R_0.
\end{equation}
Here $R_0\in\R^{m\times m}$ is the upper-triangular factor
of the thin QR factorization, and $i=1,\ldots,n$ and $j=1,\ldots,m$ index the rows and columns of $G_0$, respectively. Throughout this section,
$I_d$ is the $d\times d$ identity, $\|\cdot\|_2$ is the
Euclidean vector norm, and $\|\cdot\|_F$ is the Frobenius matrix norm.
Two-pass Gram--Schmidt orthogonalization produces
$S_0=[s_{1,0},\ldots,s_{m,0}]\in\R^{n\times m}$
with $S_0^\top S_0=I_m$, up to floating-point error.
A separate random stream prevents subspace construction from consuming
the random draws used by \texttt{PSO}.

At a successful population or \texttt{TR} step, let $x^-$ and $x^+$ be
the incumbent before and after acceptance, and let $S$ be the
current basis. Define the displacement and the column selected
for replacement by
\begin{equation}\label{eq:subspace-displacement}
	d=x^+-x^-,\qquad
	j=\min\left(\argmaxx_{1\le i\le m}|s_i^\top d|\right).
\end{equation}
Thus, ties in the maximum-correlation selection are resolved by
the first maximizing index.
The replacement is skipped, leaving $S$ unchanged, if
\[
\|d\|_2\le10^{-12}\max\{1,\|x^-\|_2\}.
\]
Otherwise, normalize $u=d/\|d\|_2$ and let $S_{-j}$ contain
all columns except $s_j$. Compute the candidate residual
$\widetilde s_j$ by reorthogonalizing $u$ against the columns
of $S_{-j}$ in two passes.
The candidate is retained only if $\|\widetilde s_j\|_2>10^{-10}$;
otherwise, $S$ remains unchanged.
For a retained candidate, the exact-arithmetic formulas are
\begin{equation}\label{eq:subspace-column}
	\widetilde s_j=(I_n-S_{-j}S_{-j}^\top)u,
	\qquad
	s_j^+=\frac{\widetilde s_j}{\|\widetilde s_j\|_2}.
\end{equation}
All other columns are preserved, and in exact arithmetic,
\begin{equation}\label{eq:subspace-preserve}
	s_i^+=s_i\quad(i\ne j),\qquad
	(S^+)^\top S^+=I_m.
\end{equation}
The computed basis satisfies the orthonormality condition up to
floating-point error.

A performed replacement preserves $m-1$ previous basis columns
and ensures that $d$ lies in the span of $S^+$ in exact arithmetic.
Population displacements with a nonzero component outside the
current subspace can therefore change its span.
For an accepted \texttt{TR} displacement $d=Sz$, with $z\in\R^m$, however,
\[
(I_n-S_{-j}S_{-j}^\top)d=z_js_j.
\]
Since the selected coefficient $z_j$ is nonzero for $d\ne0$,
a performed replacement reduces to
$s_j^+=\operatorname{sgn}(z_j)s_j$ in exact arithmetic.
Thus, success-driven updates after accepted \texttt{TR} steps can change
the sign of the selected column but do not introduce a new subspace
direction. Such updates can occur several times within one
successive-\texttt{TR} invocation, immediately after its accepted steps.
The sign changes can reverse subsequent forward-probe directions
even though the subspace is unchanged.

If the global best does not improve for five consecutive completed
generations, the last basis column is replaced by a newly sampled
Gaussian direction orthogonalized against the others.
This refresh requires no objective evaluation and is separate
from the success-driven update.
When $m=n$, the basis spans $\R^n$, and both the success-driven
updates and these one-column random refreshes reduce to possible
column-sign changes in exact arithmetic.
Neither the random initialization nor the basis-update rule alone
provides a convergence guarantee.

\subsection{Direct reduced-curvature construction and incremental updating}
\label{subsec:pso-B}

We construct an $m\times m$ curvature surrogate directly in the
current search subspace, without forming an ambient $n\times n$
curvature matrix. At the first eligible \texttt{TR} invocation, take a
snapshot of the current swarm positions $X$, velocities $V$,
and personal bests $P$, each in $\R^{N\times n}$, and form
\begin{equation}\label{eq:pso-geometry}
	D=\begin{bmatrix}
		X-\mathbf1 g^\top\\ V\\ P-X
	\end{bmatrix}\in\R^{3N\times n}.
\end{equation}
Here, $\mathbf1\in\R^N$ is the all-ones column vector,
and the rows of $X,V,P$ are $x_i^\top,v_i^\top,p_i^\top$
from the snapshot, respectively.
The vector $g\in\R^n$ is the swarm best at that snapshot,
distinct from the generation index.
Put $s=\max_{i,j}|D_{ij}|$ and $\widetilde D=D/s$ when
$s>10^{-15}$; otherwise take $\widetilde D=0$.
Let $\|\cdot\|_F$ denote the Frobenius matrix norm. For
$\widetilde D\in\R^{3N\times n}$,
\[
\|\widetilde D\|_F
=\left(\sum_{i=1}^{3N}\sum_{j=1}^{n}
|\widetilde D_{ij}|^2\right)^{1/2}.
\]
The stored geometry factor and regularization weight are
\begin{equation}\label{eq:pso-factor}
	(F,\alpha)=\begin{cases}
		\bigl(\sqrt{0.95n/\|\widetilde D\|_F^2}\,
		\widetilde D,\,0.05\bigr),
		&\widetilde D\ne0,\\
		\bigl(\mathbf0_{0\times n},\,1\bigr),
		&\widetilde D=0.
	\end{cases}
\end{equation}
Here, $\mathbf0_{0\times n}$ denotes the empty matrix with
zero rows and $n$ columns. The scalar $\alpha$ in this subsection
is a regularization weight, distinct from the socket angle
used in the contact model.
This normalization preserves direction information in the
frozen swarm snapshot, while the second case avoids a
degenerate normalization.

The operator $C_{\rm geo}=F^\top F+\alpha I_n$ is useful
for interpreting the construction, but is never formed
in the implementation. Instead, store
\begin{equation}\label{eq:pso-reducedB}
	T=FS,\qquad B=T^\top T+\alpha I_m.
\end{equation}
Because $\alpha>0$, the reduced matrix is positive definite
in exact arithmetic. It is a regularized swarm-geometry
surrogate, not an objective-Hessian estimate.
The factor $F$ and weight $\alpha$ are frozen after their
initial construction; subsequent changes to $B$ result
only from changes to $S$.

In this subsection, $t_i=T(:,i)$ denotes the $i$th factor
column, distinct from the contact tangent used in the
planar contact model.
For a nonempty factor, $T\in\R^{3N\times m}$ and
$t_i\in\R^{3N}$.
In the empty-factor case, $T\in\R^{0\times m}$;
each column has length zero, its inner products are zero,
and $B=I_m$.
The Kronecker delta satisfies $\delta_{ij}=1$ for $i=j$
and zero otherwise.

When column $j$ of $S$ is replaced, update only the
corresponding column of $T$:
\begin{equation}\label{eq:pso-Tupdate}
	t_j^+=Fs_j^+,\qquad t_i^+=t_i\quad(i\ne j).
\end{equation}
The affected row and column of the reduced curvature matrix
then satisfy
\begin{equation}\label{eq:pso-Bupdate}
	(B^+)_{ij}=(t_i^+)^\top t_j^+
	+\alpha\delta_{ij},\qquad
	(B^+)_{ji}=(B^+)_{ij},\quad 1\le i\le m.
\end{equation}
The $(m-1)\times(m-1)$ submatrix indexed by $i,\ell\ne j$
remains unchanged. The same update applies to a random-column
refresh. This exact one-row/one-column modification relies on
the frozen factor and regularization weight:
updating the swarm snapshot itself would generally change
every entry of $B$.

The matrices $F$, $S$, $T$, and $B$ require
$O(Nn+nm+m^2)$ storage, apart from constant factors.
Replacing a basis column and updating the reduced curvature
requires $O(nm+Nn+Nm)$ arithmetic.
No ambient $n\times n$ curvature matrix is formed.
In exact arithmetic, $S^\top S=I_m$ gives
$B=S^\top C_{\rm geo}S$; for the computed basis,
this relation holds up to floating-point error.
It does not assert that the frozen geometry surrogate
tracks later swarm configurations.

\subsection{Derivative-free quadratic models and successive \texttt{TR} steps}
\label{subsec:pso-trsolver}

Consider an eligible \texttt{TR} invocation after a non-improving
\texttt{PSO} generation. Index its reduced subproblems by
$k=0,1,\ldots$. Let $x_k\in\R^n$ be the current accepted center,
and let $S_k\in\R^{n\times m}$ and $B_k\in\R^{m\times m}$
be the current persistent basis and reduced curvature matrix.
Initially $x_0=g_{g+1}^{\rm PSO}$ and the working radius is
$\Delta_0>0$.
The vector $z\in\R^m$ contains reduced step coordinates,
so the full-space displacement is $S_kz$; this vector is
distinct from the physical depth in the planar insertion model.
For $j=1,\ldots,m$, $e_j\in\R^m$ is the $j$th canonical
coordinate vector.
The positive $\Delta_k$ is the working \texttt{TR} radius,
and $h_k$ is the forward-probe length defined below.

At inner step $k$, define the reduced quadratic model
\begin{equation}\label{eq:trmodel}
	q_k(z)=f(x_k)+c_k^\top z+\tfrac12 z^\top B_kz,
	\qquad \|z\|_2\le\Delta_k.
\end{equation}
Only $c_k\in\R^m$ is fitted from fresh function values;
$B_k$ is the reduced surrogate described above.
For a minimum radius $\Delta_{\min}>0$, choose
\[
h_k=\min\{\Delta_k,\max\{\Delta_k/2,\Delta_{\min}\}\}.
\]
The forward probes and curvature-corrected differences are
\begin{equation}\label{eq:tr-probes}
	f_{k,j}=f(x_k+h_kS_ke_j),\qquad
	y_{k,j}=f_{k,j}-f(x_k)-\tfrac12 h_k^2(B_k)_{jj},
	\quad j=1,\ldots,m.
\end{equation}
These $m$ probes use the current basis and the current
accepted center. Any basis-column update following an accepted
\texttt{TR} step is therefore reflected in the next probes and
fitted model.

With $A_k=h_kI_m$,
$y_k=(y_{k,1},\ldots,y_{k,m})^\top$, and a nonnegative
ridge weight $\lambda_{\rm ridge}$, the fitting problem is
\begin{equation}\label{eq:tr-fit-ridge}
	c_k=\operatorname*{arg\,min}_{c\in\R^m}
	\left\{\|A_kc-y_k\|_2^2
	+\lambda_{\rm ridge}h_k^2\|c\|_2^2\right\}.
\end{equation}
Its solution is explicit (with $c_k=\nabla_z q_k(0)$ denoting
the model gradient in reduced coordinates):
\begin{equation}\label{eq:trfit}
	(c_k)_j=
	\frac{f(x_k+h_kS_ke_j)-f(x_k)
		-\tfrac12 h_k^2(B_k)_{jj}}
	{h_k(1+\lambda_{\rm ridge})}.
\end{equation}
Each inner subproblem therefore requires at most $m$ fresh
forward probes and, when a valid descent model is obtained,
one trial evaluation. Invalid or out-of-box probes are charged
to the evaluation budget and invalidate that particular model.
The optimization is over a Euclidean ball in reduced coordinates,
not the exact intersection with the original box; the candidate
trial is checked by the objective oracle, which assigns an
invalid result to an infeasible point. No clipping of a
\texttt{TR} trial back into the box is assumed.

The reduced quadratic subproblem is
\begin{equation}\label{eq:trsub}
	z_k^*\in\operatorname*{arg\,min}_{\|z\|_2\le\Delta_k}
	\left(c_k^\top z+\tfrac12 z^\top B_kz\right).
\end{equation}
Since $B_k\succ0$ in exact arithmetic, this is a strictly
convex problem with a unique minimizer. The \texttt{TR}
optimality conditions~\cite{moresorensen1983} state that
a multiplier $\lambda_k\ge0$ exists such that
\begin{equation}\label{eq:tr-kkt}
	\begin{aligned}
		&(B_k+\lambda_kI_m)z_k^*=-c_k,\\
		&\|z_k^*\|_2\le\Delta_k,\qquad
		\lambda_k(\|z_k^*\|_2-\Delta_k)=0.
	\end{aligned}
\end{equation}
These conditions are necessary and sufficient for the reduced
quadratic problem. Writing
\begin{equation}\label{eq:tr-eig}
	B_k=U_k\operatorname{diag}(\beta_{k,1},\ldots,
	\beta_{k,m})U_k^\top,\qquad a_k=U_k^\top c_k,
\end{equation}
with $U_k^\top U_k=I_m$ and $\beta_{k,i}>0$,
the unconstrained minimizer $-B_k^{-1}c_k$ is used if it
lies inside the radius. Otherwise $\lambda_k>0$ and
\begin{equation}\label{eq:tr-spectral-solution}
	z_k^*(\lambda_k)=-U_k\operatorname{diag}\!\left(
	\frac{1}{\beta_{k,1}+\lambda_k},\ldots,
	\frac{1}{\beta_{k,m}+\lambda_k}\right)a_k.
\end{equation}
The multiplier is the unique positive root of the secular equation
\begin{equation}\label{eq:tr-secular}
	\Phi_k(\lambda)=
	\sum_{i=1}^m\frac{(a_k)_i^2}{(\beta_{k,i}+\lambda)^2}
	-\Delta_k^2=0,\qquad\lambda>0.
\end{equation}
The implementation brackets and bisects this root and applies
a final radius projection when needed. Its numerical kernel
also handles singular or indefinite input defensively, although
these cases are not expected for the regularized matrix
in~\eqref{eq:pso-reducedB}. Solving the reduced quadratic
problem does not imply global optimality for the original
BBO objective.

For a candidate $x_k^{\rm tr}=x_k+S_kz_k^*$,
the predicted decrease is
\begin{equation}\label{eq:tr-pred}
	\operatorname{pred}_k=-c_k^\top z_k^*
	-\tfrac12(z_k^*)^\top B_kz_k^*.
\end{equation}
The trial is evaluated only when the predicted decrease is
finite and positive. The actual-to-predicted reduction ratio is
\begin{equation}\label{eq:traccept}
	\rho_k=\frac{f(x_k)-f(x_k^{\rm tr})}
	{\operatorname{pred}_k}.
\end{equation}
The \texttt{TR} ratio $\rho_k$ measures agreement between
the actual objective reduction and the reduction predicted
by the model. The step is accepted if its objective value
and $\rho_k$ are finite and
\begin{equation}\label{eq:tr-accept-conditions}
	f(x_k^{\rm tr})<f(x_k)-\varepsilon_{\rm imp},
	\qquad \rho_k\ge\eta_{\rm accept},
\end{equation}
where $\varepsilon_{\rm imp}\ge0$ is the absolute improvement
tolerance and $\eta_{\rm accept}$ is the acceptance threshold.

Upon acceptance, set $x_{k+1}=x_k^{\rm tr}$, retain
$\Delta_{k+1}=\Delta_k$, and apply the success-driven
basis-column update using $x_{k+1}-x_k$, as described in
Section~\ref{subsec:pso-subspace}.
If a column changes, update the corresponding row and column
of the reduced matrix before fitting the next quadratic model.
Accepted steps can therefore be followed by further
\texttt{TR} subproblems within the same invocation.

If a subproblem is unsuccessful {\bf before any \texttt{TR}
	success} in the current invocation, the incumbent and reduced
state are retained:
\[
x_{k+1}=x_k,\qquad S_{k+1}=S_k,\qquad B_{k+1}=B_k.
\]
If the minimum radius and the remaining per-invocation and
global evaluation allowances permit another attempt, set
$\Delta_{k+1}=\Delta_k/2$ and continue; otherwise, the
invocation terminates.
Once at least one \texttt{TR} step has succeeded, the
{\bf first subsequent unsuccessful subproblem} terminates
the invocation without accepting its trial.
Exhaustion of the per-invocation or global evaluation
allowance also terminates the inner sequence.
The configured initial radius $\Delta_0$ is not modified
by these internal reductions.

If the invocation contains accepted steps, its final accepted
point is installed as the \texttt{PSO} global best.
The particle with the worst personal-best objective value has
both its position and personal best replaced by this point,
and its velocity is reset to zero. If there is no accepted \texttt{TR} step, the pre-invocation
global best is retained. The next major iteration resumes the usual \texttt{PSO} population updates while the global evaluation budget permits; there is no blending between swarm and \texttt{TR} points.

The experimental configurations and evaluation allowances are
reported in Section~\ref{sec:bcp-profiles} and
Table~\ref{tab:psotr-experimental-settings}.

\section{\texttt{PSO+TR} robotic experimental design, controller selection, and validation}
\label{sec:robotic-selection}

To test whether sampled training gains transfer to new scenarios, we separate optimization, candidate screening, frozen-policy validation, and withheld final testing. The
proposed \texttt{PSO+TR} solver optimizes a prescribed 14-parameter
controller, while a separate screen--freeze--validate protocol determines
whether a tuned policy should replace the fixed staged reference.
This separation is essential: an improvement in the sampled training
objective does not, by itself, establish reliable performance on
independent scenarios. The complete procedure comprises optimization,
candidate screening, frozen-policy validation, and withheld final tests.

\subsection{Reconstructed simulator and optimization protocol}
\label{subsec:robotic-experimental-protocol}

The BBO oracle is implemented in \texttt{MATLAB} using the four-corner,
nearest-face contact model and the controller defined in Section \ref{sec:PlanarInsertion}. The simulator uses semi-implicit Euler integration with a
$0.250$~ms physics step for training and a $0.125$~ms step for
screening, validation, and final testing. Control commands are updated
every $10$~ms, and each episode has a maximum duration of $8$~s.
This is a synthetic numerical implementation of the stated planar model;
it is not a hardware experiment or a verified reproduction of an
independently implemented simulator.

Eight independent optimizer runs use seeds $0,\ldots,7$, ten swarm
particles, a 14-dimensional normalized parameter vector, and a budget
of $120$ charged objective-evaluation attempts per seed.
Each feasible evaluation executes eight nominal and eight development
scenarios, which are drawn once for a given seed and held fixed
throughout that optimizer run. The training oracle is the maximum
of the two sample means of the penalized episode cost
in~\eqref{eq:problem}.

The charged training workload is
$8\times120\times16=15{,}360$ rollout equivalents,
or $1{,}920$ per seed. An infeasible attempted evaluation is charged
by the optimizer without executing its scenario batch.
The saved optimizer logs record $109$ such infeasible attempts
across the eight runs, leaving $960-109=851$ feasible evaluations
and $851\times16=13{,}616$ executed training episodes.
Thus, charged rollout equivalents and executed simulator episodes
are reported separately.

Initialization evaluations, subsequent particle evaluations,
\texttt{TR} model probes, and \texttt{TR} trial evaluations all count against
the same $120$-attempt budget, including infeasible attempts.
The successive-\texttt{TR} mechanism described in Section~\ref{sec:pso-tr}
may solve more than one reduced subproblem within a conditional
\texttt{TR} invocation; it receives no additional unaccounted
objective-evaluation allowance.

The \texttt{PSO+TR} hyperparameters are transferred from the separate
\texttt{BCP} tuning archive without selecting or readjusting them using the
robotic final-test results. The controller parameters themselves are
optimized independently for each seed. The \texttt{BCP} solver comparison is
reported separately in Section~\ref{sec:bcp-profiles}; it is not a
matched-simulator robotic comparison among optimizers.

\subsection{Candidate screening and policy freezing}
\label{subsec:robotic-screening}

Following optimization within each seed, the returned incumbent is
placed first in the candidate list, followed by archived population
points with finite training scores in increasing order of
$\widehat f^{(s)}$ from~\eqref{eq:pso-score}.
The archive contains finite evaluations of the initial and subsequent
swarm positions, recorded in evaluation order. Equal training scores
preserve this archive order. \texttt{TR} probes and trials are not appended
to the population archive; the returned incumbent is included
separately even if it results from an accepted \texttt{TR} step.
Exact duplicate parameter vectors are removed from the combined list
while preserving their first occurrence. The first $13$ distinct
vectors are screened, or all available vectors if fewer remain.

Within each seed, every collected candidate and the fixed staged
reference are evaluated on the same newly sampled $64$ nominal
and $64$ development scenarios. These scenarios are independent
of the optimizer's fixed training batch.

For a candidate $\theta$, let $\widehat v_j(\theta)$ be its empirical
violation fraction on screening distribution $j\in\{0,1\}$, and let
$\widehat\ell_j(\theta)$ be the fraction of paired screening scenarios
on which the fixed staged controller succeeds but $\pi_\theta$ fails.
A candidate is eligible only if
\begin{equation}\label{eq:robotic-screening}
	\widehat v_j(\theta)\leq0.05,
	\qquad
	\widehat\ell_j(\theta)\leq0.03,
	\qquad j\in\{0,1\}.
\end{equation}
Among eligible candidates, the method selects the one with the lowest
balanced empirical mean episode cost,
$\tfrac12(\overline C_0(\theta)+\overline C_1(\theta))$, where
\[
\overline C_j(\theta)=\frac{1}{64}\sum_{r=1}^{64}
C(\pi_\theta,\omega^{\rm scr}_{j,r}),
\qquad j\in\{0,1\}.
\]
Here $\omega^{\rm scr}_{j,r}$ is the $r$th screening scenario
from distribution $P_j$ within the current seed.
If eligible candidates have exactly equal balanced costs, the first
one in the candidate-list order is retained; no tie tolerance is used.

This selection criterion excludes the additional penalties
$10(1-S)+20R$ used in the training oracle. The episode cost $C$
retains its own unsuccessful-completion and protective-abort penalties.
The screening inequalities separately restrict the observed violation
and paired-loss frequencies. The selected candidate is then frozen;
no validation observation is used to change its parameters or to
reselect a different candidate. If no candidate is eligible, the
fixed staged reference is retained without candidate validation.

\subsection{Independent frozen-policy validation}
\label{subsec:robotic-frozen-validation}

The frozen candidate and the fixed staged reference are evaluated on a
new, common batch of $512$ nominal and $512$ development scenarios.
For $k$ occurrences of a binary event in $n$ independent scenarios
drawn from a given validation distribution, we use the one-sided
exact Clopper--Pearson upper confidence bound
\begin{equation}\label{eq:robotic-cp}
	U(k,n,a_\star)=
	\begin{cases}
		\operatorname{Beta}^{-1}(1-a_\star;k+1,n-k),&k<n,\\
		1,&k=n.
	\end{cases}
\end{equation}
Here, $\operatorname{Beta}^{-1}(q;a,b)$ is the $q$ quantile
of the beta distribution with shape parameters $a,b>0$.
The fixed per-bound tail probability is $a_\star=0.05/256$.
Each upper bound therefore has nominal confidence level
$1-a_\star$, with coverage at least this level under the
binomial model.
In validation, $n=512$ and $k$ counts occurrences of the
event being bounded within the relevant distribution.
The denominator $256$ is a prescribed allocation factor;
it is not presented as the number of simultaneous tests
in the eight-seed experiment.

For {\bf each} validation distribution, adoption requires the bound
for the paired event ``fixed staged succeeds, candidate fails'' to be
at most $0.03$, the candidate's violation bound to be at most $0.05$,
and the fixed staged reference's violation bound to be at most $0.05$.
In addition, the candidate's balanced mean validation cost must be
strictly below that of the reference. All these conditions must hold
for the {\bf same frozen candidate}. If any condition fails, or no
candidate passes screening, the deployed policy is the fixed staged
reference. Rejected candidates are never replaced by the zero-parameter
hybrid controller, because enabling its probing and retreat rules does
not reproduce the fixed staged baseline.

The exact binomial bound quantifies uncertainty about the specified
binary event probabilities under independent draws from the validation
distributions. It does not give a guarantee for the withheld stress
distribution, other contact models, or physical hardware; the empirical
cost-improvement condition is not itself a confidence-bound guarantee.

\subsection{Final tests, evaluation accounting, and provenance}
\label{subsec:robotic-final-tests}

After deployment is determined, four controller configurations are
evaluated: the deployed \texttt{PSO+TR} selection pipeline (the adopted
candidate or its fixed staged fallback), fixed staged, continued search, and fixed hybrid.
The continued-search reference uses the nominal staged parameters
and disables probing and retreat, but continues the lateral sinusoid
during insertion instead of setting its lateral target to the measured
current position; its insertion angular target remains the measured angle. Each configuration is evaluated on the {\bf same}
$400$ fresh scenarios per seed for each of the nominal, development,
and withheld stress distributions. Scenarios are paired across
controllers within a seed, and fresh draws are used for different
seeds. The resulting final-test ledger contains
$8\times4\times3\times400=38{,}400$ controller episodes, not
$38{,}400$ independent optimizer runs.

The archive includes the seed-specific optimizer checkpoints,
function-evaluation histories, screening and validation decisions,
and episode-level final-test outcomes. The supplied results also contain an auditable correction of an earlier deployment error. Seed~4 produced no eligible screening candidate, so the prescribed protocol required the fixed staged fallback without candidate validation. The original runner instead
evaluated the zero-parameter hybrid for this seed. Its reported
deployed-policy final-test rows were deterministically replaced by
the {\bf already measured} fixed staged rows on exactly matched
scenarios. The original records and a fallback-correction ledger
are retained; this correction does not represent additional simulations.

The complete offline episode budget and final-test count are summarized
in Table~\ref{tab:psotr-native-budget}; the table is regenerated from
the experimental manifest and decision records:


\begin{table}[!htbp]\centering\small
	\caption{Native \texttt{PSO+TR} reconstructed-simulator episode
		accounting across all eight optimizer seeds. Executed training
		episodes include feasible \texttt{TR} probes and trials.
		The charged training workload is 15,360 rollout equivalents;
		infeasible attempts do not execute episodes. Validation is
		skipped for seed~4, which has no eligible screening candidate.
		The final-test count refers to corrected archive records.}
	\label{tab:psotr-native-budget}
	\begin{tabularx}{\linewidth}{Xr}\toprule
		Stage & Count \\
		\midrule
		Training: 851 executed objective calls $\times$ 16 episodes
		& 13616 \\
		Candidate screening: 104 candidates $\times$ 2 groups $\times$ 64
		& 13312 \\
		Fixed staged screening: 8 seeds $\times$ 2 groups $\times$ 64
		& 1024 \\
		Frozen-candidate validation: 7 seeds $\times$ 2 groups $\times$ 512
		& 7168 \\
		Fixed staged validation: 7 seeds $\times$ 2 groups $\times$ 512
		& 7168 \\
		\midrule
		Offline total: executed episodes
		& 42288 \\
		\midrule
		Final-test records: 8 seeds $\times$ 4 policies
		$\times$ 3 groups $\times$ 400
		& 38400 \\
		\bottomrule
	\end{tabularx}
\end{table}

The numerical results, controller comparisons, and \texttt{TR} diagnostics are
reported next. An adopted controller is a policy that satisfies the
independent deployment criteria, whereas an accepted \texttt{TR} step is an
internal optimization event; the two statistics measure different
outcomes. Demonstrating the incremental benefit of the \texttt{TR} safeguard
would require a matched robotic \texttt{PSO}-only ablation with the
same simulation model, evaluation budget, scenario design, and
selection protocol.


\section{\texttt{PSO+TR} robotic controller-tuning results}
\label{sec:psotr-robot-results}

To assess the complete tuning and deployment pipeline, we compare
its final-test outcomes with three fixed controllers on paired
nominal, development, and withheld stress scenarios.

Eight independent optimization runs tune the same 14-parameter hybrid
controller on the reconstructed \texttt{MATLAB} planar contact
simulator. The historical robotic \texttt{MAES} measurements
are not pooled with these \texttt{PSO+TR} results.
Each run uses ten particles, the independently selected \texttt{TR}
configuration, and a fixed training batch of eight nominal
and eight development scenarios shared across feasible
candidate evaluations within that run.

The optimization budget is $120$ charged evaluation attempts
per seed, or $1{,}920$ charged training rollout equivalents.
Across the eight runs, $109$ of the $960$ charged attempts are
infeasible and do not execute the simulator batch. The remaining
$851$ feasible evaluations execute $13{,}616$ training episodes.
Initialization evaluations, particle evaluations, \texttt{TR} model probes,
and \texttt{TR} trials share this budget, including infeasible attempts.

Each seed screens $13$ distinct candidates selected from the returned
incumbent and the finite population archive. All candidates and the
fixed staged reference use the same $64$ fresh scenarios per nominal
and development distribution within that seed. When screening selects
an eligible candidate, that frozen candidate and the reference are
validated on $512$ new paired scenarios per distribution.
Final testing uses $400$ new scenarios per seed, distribution,
and controller. The final-test archive contains $38{,}400$ records
across four controllers. Scenario streams are independent across
seeds, while controller comparisons within each seed use identical
test scenarios.

Seed~4 produced no eligible screening candidate and therefore
used the fixed staged fallback without candidate validation.
The original run evaluated the zero-parameter hybrid
controller as the fallback. The published dataset replaces the
corresponding \texttt{PSO+TR} records with the already measured
fixed staged outcomes on the identical paired final-test scenarios.
The unmodified seed checkpoints and deterministic correction audit
are retained. The corrected \texttt{MATLAB} runner explicitly
implements the fixed staged fallback.

\begin{figure}[!htbp]
	\centering
	\includegraphics[width=.96\linewidth]{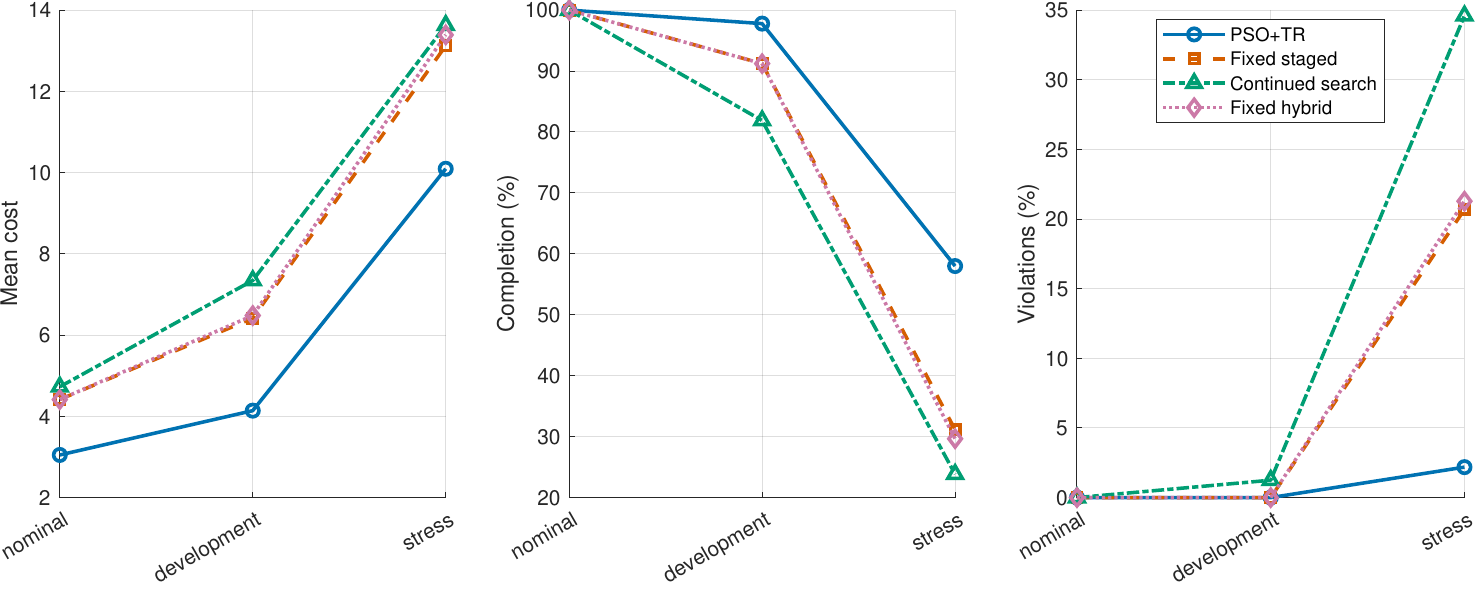}
	\caption{Deployed-policy outcomes on paired final-test scenarios.}
	\label{fig:psotr-robot-performance}
\end{figure}

Table~\ref{tab:psotr-native-robotic} reports final-test outcomes
as means and sample standard deviations across the eight optimizer
seeds. Figure~\ref{fig:psotr-robot-performance} displays the
corresponding controller comparisons. Cost averages include
successful, unsuccessful, and fallback outcomes.

\begin{table}[htbp]
	\centering\small
	\caption{Robotic final-test outcomes across eight optimizer seeds.}
	\label{tab:psotr-native-robotic}
	\begin{tabular}{llrrr}
		\toprule
		Group & Method & Completion (\%) & Violation (\%) & Cost \\
		\midrule
		nominal & \texttt{PSO+TR}
		& 100.00 $\pm$ 0.00 & 0.00 $\pm$ 0.00 & 3.050 $\pm$ 0.581 \\
		nominal & Fixed staged
		& 100.00 $\pm$ 0.00 & 0.00 $\pm$ 0.00 & 4.414 $\pm$ 0.032 \\
		nominal & Continued search
		& 100.00 $\pm$ 0.00 & 0.00 $\pm$ 0.00 & 4.728 $\pm$ 0.036 \\
		nominal & Fixed hybrid
		& 100.00 $\pm$ 0.00 & 0.00 $\pm$ 0.00 & 4.414 $\pm$ 0.032 \\
		development & \texttt{PSO+TR}
		& 97.78 $\pm$ 3.52 & 0.00 $\pm$ 0.00 & 4.140 $\pm$ 1.033 \\
		development & Fixed staged
		& 91.22 $\pm$ 1.22 & 0.00 $\pm$ 0.00 & 6.415 $\pm$ 0.176 \\
		development & Continued search
		& 81.88 $\pm$ 2.19 & 1.25 $\pm$ 0.65 & 7.345 $\pm$ 0.213 \\
		development & Fixed hybrid
		& 91.22 $\pm$ 1.22 & 0.00 $\pm$ 0.00 & 6.482 $\pm$ 0.177 \\
		stress & \texttt{PSO+TR}
		& 58.00 $\pm$ 13.05 & 2.19 $\pm$ 5.99 & 10.094 $\pm$ 1.789 \\
		stress & Fixed staged
		& 31.12 $\pm$ 1.87 & 20.69 $\pm$ 2.54 & 13.124 $\pm$ 0.312 \\
		stress & Continued search
		& 23.84 $\pm$ 1.89 & 34.59 $\pm$ 2.22 & 13.628 $\pm$ 0.302 \\
		stress & Fixed hybrid
		& 29.66 $\pm$ 2.20 & 21.28 $\pm$ 2.61 & 13.391 $\pm$ 0.333 \\
		\bottomrule
	\end{tabular}
\end{table}

Under development conditions, the deployed \texttt{PSO+TR}
pipeline achieves $97.78\%$ completion versus $91.22\%$ for
fixed staged. On the stress scenarios, completion is $58.00\%$
versus $31.12\%$, a paired difference of $26.88$ percentage points.
The corresponding stress violation fractions are $2.19\%$ and
$20.69\%$, and mean costs are $10.094$ and $13.124$, respectively.
These are descriptive contrasts across eight optimization seeds,
not simultaneous confidence bounds or a direct \texttt{PSO}-only
robotic ablation.

\begin{figure}[!htbp]
	\centering
	\includegraphics[width=.96\linewidth]{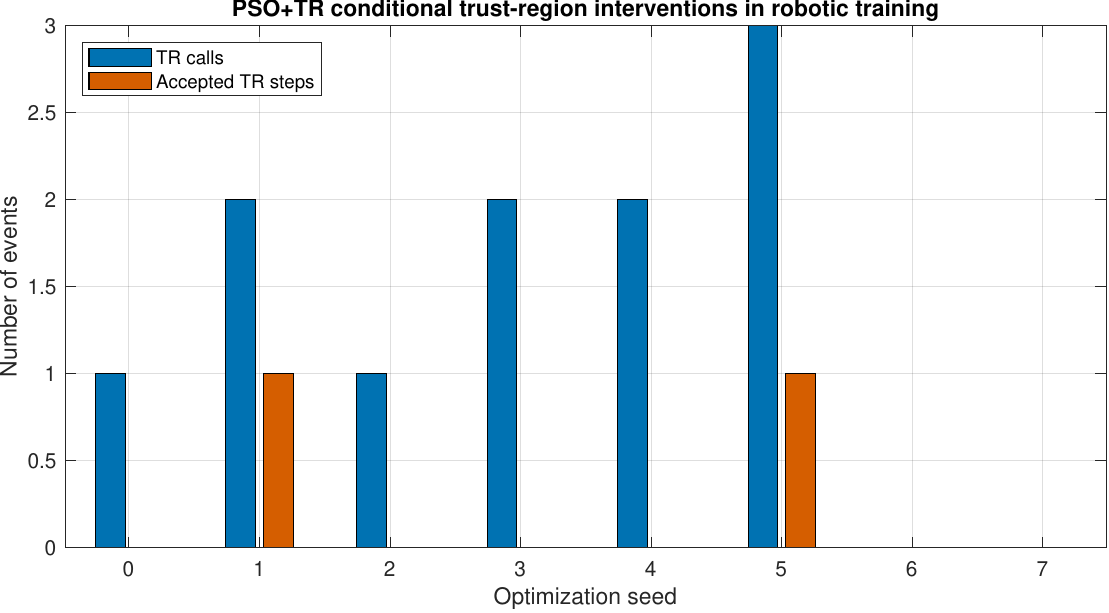}
	\caption{Conditional \texttt{TR} invocations and accepted
		steps by optimizer seed.}
	\label{fig:psotr-robot-tr-diagnostics}
\end{figure}

Table~\ref{tab:psotr-native-decisions} lists the seed-level
optimization and deployment decisions, while
Figure~\ref{fig:psotr-robot-tr-diagnostics} displays the
\texttt{TR} intervention counts.
Seven of the eight optimization runs produced a screening-eligible
candidate. All seven frozen candidates were adopted after independent
validation; seed~4 used the fixed staged fallback.
Across all seeds, \texttt{PSO+TR} attempted $11$ conditional \texttt{TR}
invocations and accepted $2$ \texttt{TR} steps, consuming $226$ charged
\texttt{TR} evaluation attempts. The accepted steps occurred in
seeds~1 and~5.
Adoption of a final controller and acceptance of a \texttt{TR}
step are separate events. These sparse \texttt{TR} improvements
do not identify a causal benefit of \texttt{TR} relative to a
matched \texttt{PSO}-only robotic experiment.

\begin{table}[!htbp]
	\centering\small
	\caption{Seed-level optimization and controller-adoption decisions.}
	\label{tab:psotr-native-decisions}
	\begin{tabular}{rrrrrr}
		\toprule
		Seed & Adopted & Charged evals & \texttt{TR} calls
		& \texttt{TR} accepted & Charged \texttt{TR} evals \\
		\midrule
		0 & 1 & 120 & 1 & 0 & 20 \\
		1 & 1 & 120 & 2 & 1 & 50 \\
		2 & 1 & 120 & 1 & 0 & 10 \\
		3 & 1 & 120 & 2 & 0 & 40 \\
		4 & 0 & 120 & 2 & 0 & 40 \\
		5 & 1 & 120 & 3 & 1 & 66 \\
		6 & 1 & 120 & 0 & 0 & 0 \\
		7 & 1 & 120 & 0 & 0 & 0 \\
		\midrule
		Total & 7 & 960 & 11 & 2 & 226 \\
		\bottomrule
	\end{tabular}
\end{table}

Table~\ref{tab:psotr-native-paired} reports differences between the
deployed pipeline and fixed staged, formed within each seed before
computing their means and sample standard deviations across seeds.
Here, pp denotes percentage points; positive completion and negative
violation or cost differences favor \texttt{PSO+TR}.
Nominal completion is unchanged, while mean cost differences favor
the deployed pipeline in all three distributions.

\begin{table}[!htbp]
	\centering\small
	\caption{Paired final-test differences: deployed
		\texttt{PSO+TR} minus fixed staged.}
	\label{tab:psotr-native-paired}
	\begin{tabular}{lrrr}
		\toprule
		Distribution & $\Delta$ completion (pp)
		& $\Delta$ violation (pp) & $\Delta$ cost \\
		\midrule
		nominal
		& 0.00 $\pm$ 0.00 & 0.00 $\pm$ 0.00 & -1.365 $\pm$ 0.587 \\
		development
		& 6.56 $\pm$ 2.97 & 0.00 $\pm$ 0.00 & -2.275 $\pm$ 0.965 \\
		stress
		& 26.88 $\pm$ 12.04 & -18.50 $\pm$ 7.73 & -3.030 $\pm$ 1.584 \\
		\bottomrule
	\end{tabular}
\end{table}

\begin{figure}[!htbp]
	\centering
	\includegraphics[width=.96\linewidth]{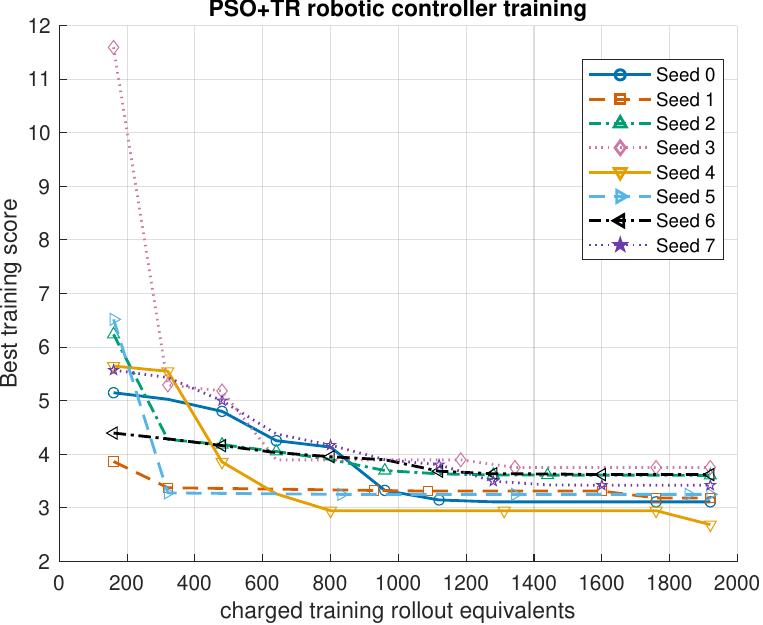}
	\caption{\texttt{PSO+TR} training scores versus charged
		training rollout equivalents.}
	\label{fig:psotr-robot-training}
\end{figure}

\begin{figure}[!htbp]
	\centering
	\includegraphics[width=.96\linewidth]{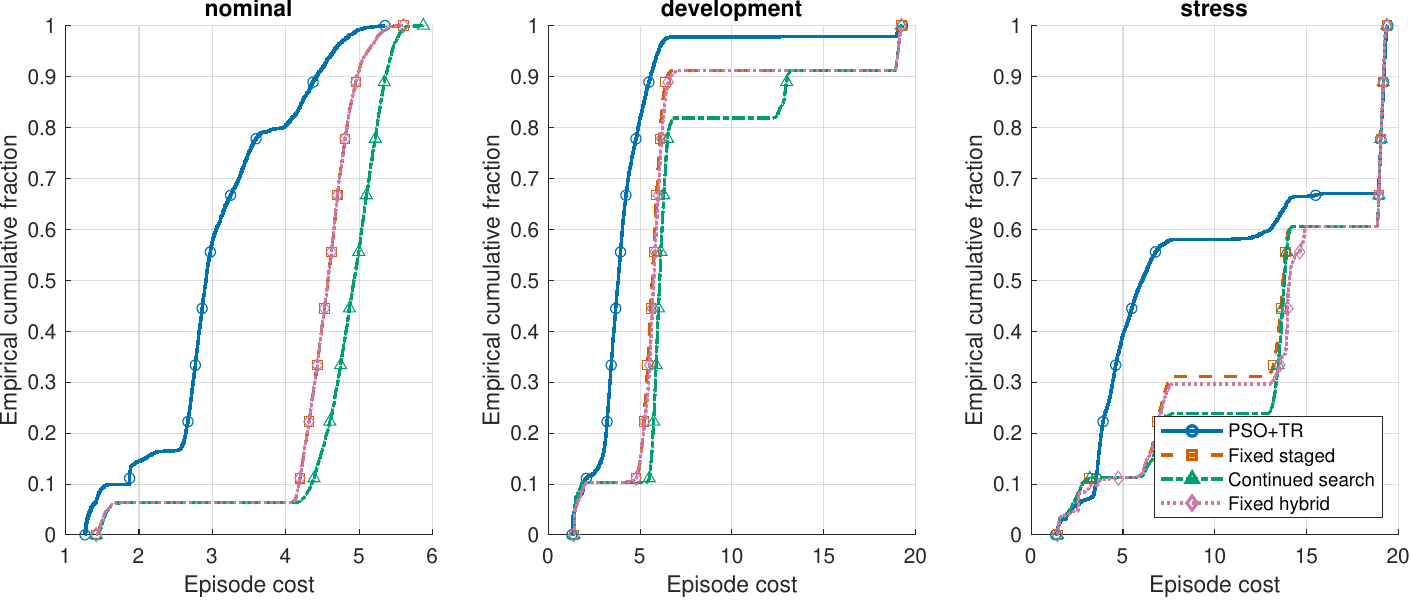}
	\caption{Empirical final-test cost distributions by controller
		and scenario group.}
	\label{fig:psotr-robot-ecdf}
\end{figure}

Figure~\ref{fig:psotr-robot-training} shows training incumbent
scores versus charged training rollout equivalents, computed as
16 times the number of charged evaluation attempts.
Infeasible attempts advance this count without executing the
simulator. Each seed uses its own fixed training batch, so the
scores refer to different sampled objectives across seeds.
Figure~\ref{fig:psotr-robot-ecdf} complements the mean-cost
comparisons with empirical episode-cost distributions by controller
and scenario group, using all $38{,}400$ final-test records,
including unsuccessful and fixed staged fallback outcomes.

The simulator is synthetic and represents local planar motion:
it omits articulated joints, gravity, sensor and actuator dynamics,
and three-dimensional contact. The tuned controller is not certified
on physical hardware. Fixed training batches support consistent
comparisons of \texttt{TR} center, probe, and trial values within each seed;
the stress distribution is reserved for final evaluation.
The results assess the complete tuning and deployment pipeline
on this simulator. They do not isolate the contribution of \texttt{TR}
or establish robustness to other contact models.

\FloatBarrier

\section{Numerical optimizer comparison on \texttt{BCP} problems}
\label{sec:bcp-profiles}

To assess target attainment and evaluation efficiency, we compare
\texttt{PSO+TR} with four other optimization methods on bound-constrained
\texttt{BCP} problems under matched objective-evaluation budgets.
We specify the benchmark eligibility criteria, solver configurations,
evaluation accounting, and empirical target definitions used in
the comparison. This experiment assesses numerical optimization
performance separately from robotic controller selection and validation.

\paragraph{Benchmark and comparison protocol.}
The \texttt{BCP} collection contains 326 problem definitions.
We evaluate 93 eligible problems with dimensions $1\le n\le12$.
These problems have continuous decision variables, finite
componentwise bounds, no additional nonlinear constraints, and
a finite, real objective value at the repaired, box-projected
initial point.
Five methods are compared: \texttt{PSO+TR},
\texttt{PSO}~\cite{pso},
\texttt{VRBBO}~\cite{kimiaei2021vrbbo},
\texttt{MAES}~\cite{maes}, and
\texttt{DE}~\cite{stornprice1997}.

Each method is tested with three independent random seeds on every
problem, yielding $93\times3=279$ matched problem--seed instances
per method and $5\times279=1{,}395$ solver runs in total. For a
given instance, all methods receive the same prescribed initial candidate
and the same maximum number of objective evaluations,
\begin{equation}\label{eq:bcp-budget}
	N_{\max}(n)=\max\{100,\lceil100n\rceil\}.
\end{equation}
Every solver evaluation attempt is charged to this allowance, including
initialization, model probes, infeasible attempts, and unsuccessful
trials of the conditional \texttt{TR} safeguard. An infeasible attempt
can be charged without calling the underlying objective. The initial
eligibility check is performed before the solver runs and is outside
their evaluation budgets. The common budget does not imply identical
internal initialization or identical evaluation overhead across solvers.

\paragraph{Evaluation accounting and solver configurations.}
Particle evaluations, forward model probes, and all \texttt{TR} trial
evaluations share the same global evaluation budget. A \texttt{TR}
invocation is eligible only after a complete non-improving generation
and is limited to 32 charged evaluation attempts within that budget.
It may consume fewer attempts when a model fails early or either
evaluation allowance is exhausted.

Table~\ref{tab:psotr-experimental-settings} collects the shared
\texttt{PSO+TR} configuration and fixed implementation constants.
The archived parameter-selection record reports a budget factor of
1,000; its configuration was transferred without retuning to this
\texttt{BCP} comparison, which uses factor 100, and to the reconstructed
\texttt{MATLAB} controller experiment. Not every constant in the table
was varied in the parameter-selection study. The present comparison
therefore measures performance under its stated budget rather than
parameter optimality for that budget. Robotic controller selection
and outcomes are reported in Sections~\ref{sec:robotic-selection}
and~\ref{sec:psotr-robot-results}. The combined stochastic method is
assessed empirically; no global or stationary-point convergence
guarantee is asserted.

\begin{table}[!htbp]\centering\small
	\caption{Shared \texttt{PSO+TR} settings in the reported \texttt{BCP}
		and robotic experiments. The first group records the transferred
		configuration; the second lists fixed implementation constants.}
	\label{tab:psotr-experimental-settings}
	\begin{tabularx}{\linewidth}{@{}Xl@{}}\toprule
		Setting & Value or rule\\
		\midrule
		\multicolumn{2}{@{}l}{{\bf Transferred configuration}}\\
		Inertia weight $\omega$ & $0.7$\\
		Cognitive and social coefficients $(c_1,c_2)$ & $(1.4,1.4)$\\
		Reduced dimension $m$ & $\min\{n,10\}$\\
		Initial radius $\Delta_0$ & $0.2$\\
		Refresh interval & Five non-improving generations\\
		Acceptance threshold $\eta_{\rm accept}$ & $0.2$\\
		Ridge weight $\lambda_{\rm ridge}$ & $10^{-12}$\\
		Per-invocation evaluation cap & $32$ charged attempts\\
		\midrule
		\multicolumn{2}{@{}l}{{\bf Fixed implementation constants}}\\
		Minimum radius $\Delta_{\min}$ & $10^{-8}$\\
		Absolute improvement tolerance $\varepsilon_{\rm imp}$ & $0$\\
		Radius multiplier on a permitted retry & $1/2$\\
		Initial velocity scale relative to box width & $0.1$\\
		Forward-probe length $h_k$ & $\min\{\Delta_k,\max\{\Delta_k/2,\Delta_{\min}\}\}$\\
		Geometry-factor normalization & $\|F\|_F^2=0.95n$ for nonempty $F$\\
		Regularization weight $\alpha$ & $0.05$; $1$ when $F$ is empty\\
		Zero-geometry cutoff $\max_{i,j}|D_{ij}|$ & $10^{-15}$\\
		Displacement rejection threshold & $10^{-12}\max\{1,\|x^-\|_2\}$\\
		Minimum retained residual norm & Strictly greater than $10^{-10}$\\
		Reorthogonalization passes & $2$\\
		Basis orthogonality-check tolerance & $10^{-8}$ in Frobenius norm\\
		Eigenvector orthogonality-check tolerance & $10^{-6}$ in Frobenius norm\\
		Numerical eigenvalue threshold & $10^{-12}\max\{1,\max_i|\beta_i|,\|a\|_2\}$\\
		Secular-root bracketing and bisection limits & $200$ and $100$ iterations\\
		Maximum Gaussian draws per random refresh & $10$\\
		\bottomrule
	\end{tabularx}
\end{table}

In the numerical eigenvalue-threshold row, $\beta_i$ are the
eigenvalues of the reduced curvature matrix, and $a=U^\top c$
contains the model's linear coefficients in its orthonormal
eigenbasis $U$.

In the \texttt{BCP} runs, \texttt{PSO+TR}, \texttt{PSO}, and
\texttt{DE} use $N=\max\{12,\min\{40,5n\}\}$ population members.
Their initial positions are sampled uniformly in the box, with the
first position replaced by the prescribed initial candidate $x^0$.
All initialization evaluations are charged. The robotic runs instead
use $N=10$, $x^0=\mathbf0_{14}$, and 120 charged attempts per seed;
the shared settings in Table~\ref{tab:psotr-experimental-settings}
are unchanged.

The \texttt{BCP} random seeds are $1,2,3$. \texttt{PSO+TR} uses
separate Mersenne Twister streams for swarm and subspace draws.
For solver seed $s$, the subspace stream uses
$(s+1{,}000{,}003)\bmod(2^{32}-1)$; the robotic runner uses solver
seed $5000+s$ for experiment labels $s=0,\ldots,7$.
Maximum-correlation ties in basis-column selection use the first
maximizing index. If the initial Gaussian basis requires rank
completion, at most $20\max\{n,m\}$ additional Gaussian draws are
attempted before canonical coordinate directions are used.

Table~\ref{tab:bcp-comparator-settings} records the principal
settings of the four comparator implementations. The separate
\texttt{PSO} comparator starts with zero velocities and retains
velocity components when positions are projected onto the box.
\texttt{PSO+TR} uses randomly initialized velocities and zeros
the affected velocity components after projection. Thus the
comparison includes implementation differences beyond the
presence of the \texttt{TR} safeguard.

\begin{table}[!htbp]\centering\small
	\caption{Comparator configurations in the supplied \texttt{BCP}
		runner. Each method uses the common budget in~\eqref{eq:bcp-budget}.}
	\label{tab:bcp-comparator-settings}
	\begin{tabularx}{\linewidth}{@{}lX@{}}\toprule
		Method & Configuration\\
		\midrule
		\texttt{PSO} & $(\omega,c_1,c_2)=(0.7,1.4,1.4)$;
		zero initial velocities; componentwise position projection.\\
		\texttt{DE} & \texttt{DE/rand/1} mutation with scale $0.7$;
		binomial crossover probability $0.85$, with one forced coordinate;
		componentwise trial projection; strict-improvement replacement.\\
		\texttt{MAES} & Offspring count $\lambda=\max\{8,4+\lfloor3\log n\rfloor\}$;
		parent count $\mu=\lfloor\lambda/2\rfloor$; initial mean $x^0$,
		transformation $I_n$, and zero evolution path;
		initial scale $0.15\operatorname{median}_j(u_j-\ell_j)$,
		with floor $10^{-12}$; projected offspring and recombined mean;
		no \texttt{TR} rescue.\\
		\texttt{VRBBO} & Archived implementation with its default tuning
		structure (\texttt{alg=5}); initial point $x^0$; evaluation limit
		$N_{\max}$; no finite time limit or objective target.
		The counted oracle enforces the box by penalizing invalid attempts.\\
		\bottomrule
	\end{tabularx}
\end{table}

For the \texttt{MAES} implementation, the normalized recombination
weights are
\[
w_j=\frac{\log(\mu+1/2)-\log j}
{\sum_{\ell=1}^{\mu}[\log(\mu+1/2)-\log\ell]},
\qquad j=1,\ldots,\mu,
\qquad \mu_{\rm eff}=\frac{1}{\sum_{j=1}^{\mu}w_j^2}.
\]
Its path, rank-one, rank-$\mu$, and step-scale coefficients are
\[
\begin{aligned}
	c_\sigma&=\frac{\mu_{\rm eff}+2}{n+\mu_{\rm eff}+5},
	&c_{\rm one}&=\frac{2}{(n+1.3)^2+\mu_{\rm eff}},\\
	c_\mu&=\min\!\left\{1-c_{\rm one},
	\frac{2(\mu_{\rm eff}-2+1/\mu_{\rm eff})}
	{(n+2)^2+\mu_{\rm eff}}\right\},
	&d_\sigma&=1+c_\sigma+2\max\!\left\{0,
	\sqrt{\frac{\mu_{\rm eff}-1}{n+1}}-1\right\}.
\end{aligned}
\]
The step-scale update uses
$\chi_n=\sqrt{2}\,\Gamma((n+1)/2)/\Gamma(n/2)$, where $\Gamma$
is the gamma function. These dimension-dependent coefficients
follow the recorded implementation and were not retuned for this
comparison.

\paragraph{Target attainment and evaluation-count profiles.}
Let $\mathcal P$ be the set of $279$ matched problem--seed instances
and $\mathcal A$ the five solvers. Here $p\in\mathcal P$ indexes an
instance and $s\in\mathcal A$ indexes a solver, rather than a random seed.
For attempted evaluation $e$, let $x_{p,s}^{(e)}$ be its trial point
and $\widehat f_{p,s}^{(e)}$ the recorded oracle value, which equals
$f_p(x_{p,s}^{(e)})$ for a valid finite evaluation. Define the
recorded best-value trace
\[
F_{p,s}(r)=\min_{1\le e\le r} \widehat f_{p,s}^{(e)},\qquad
f_{p,s}^{\rm final}=F_{p,s}(N_{p,s}^{\rm used}),
\qquad f_p^\star=\min_{s\in\mathcal A}f_{p,s}^{\rm final}.
\]
Here, $N_{p,s}^{\rm used}$ is the charged evaluation count at termination,
distinct from the first-hit count $N_{p,s}$. Invalid attempts receive
$+\infty$ in \texttt{PSO+TR} and the finite penalty
$\texttt{realmax}/10^6$ in the four comparator implementations.
The comparator box checks allow a tolerance of $10^{-10}$;
\texttt{PSO+TR} uses strict box checks.
The initial value is $f_p^0=f_p(x_p^0)$.
The trace includes finite model probes and rejected \texttt{TR} trials, so its best
value can differ from the accepted incumbent returned by the solver.
In Eq.~\eqref{eq:bcp-profiles-def}, $\#$ denotes set cardinality.
We use the instance-dependent target
\begin{equation}\label{eq:bcp-target}
	t_p=f_p^\star+0.01\max\bigl\{f_p^0-f_p^\star,
	\epsilon_{\mathrm{mach}}\max(1,|f_p^\star|)\bigr\}.
\end{equation}
Thus an algorithm reaches the target when at most 1\% of the observed
initial-to-best gap remains, up to the indicated numerical tolerance. Here, $\epsilon_{\rm mach}=2^{-52}\simeq2.220446049250313\times10^{-16}$
is the spacing at one in IEEE binary64 arithmetic, as implemented by
\texttt{eps(1)}. These targets are empirical references rather than certified global optima.

Let $N_{p,s}$ denote the number of charged evaluation attempts needed
by solver $s$ to reach $t_p$ for the first time, with $N_{p,s}=+\infty$
when the target is not attained. Define
$H_p=\min_{a\in\mathcal A}N_{p,a}$. Following Dolan and
Mor\'e~\cite{dolanmore2002} and Mor\'e and Wild~\cite{morewild2009},
the performance and data profiles are
\begin{equation}\label{eq:bcp-profiles-def}
	\begin{aligned}
		\rho_s(\tau)
		&=\frac{1}{279}\#\left\{p:H_p<+\infty,\;
		N_{p,s}\le \tau H_p\right\},\\
		\delta_s(\kappa)
		&=\frac{1}{279}\#\left\{p:
		N_{p,s}\le\kappa(n_p+1)\right\}.
	\end{aligned}
\end{equation}
Here, $\tau\ge1$, $\kappa\ge0$, and $n_p$ is the problem dimension.
Instances without target attainment by any method contribute zero to
both profiles while remaining in the common denominator. The explicit
condition $H_p<+\infty$ avoids an undefined ratio for such instances. Final solution quality is further summarized by the
instance-normalized regret
\begin{equation}\label{eq:bcp-regret}
	r_{p,s}=\frac{\max\{0,f_{p,s}^{\mathrm{final}}-f_p^\star\}}
	{\max\{1,|f_p^\star|\}}.
\end{equation}
The regret normalization allows comparison across problems with different
objective scales but can still produce large values on a few instances.

\begin{figure}[!htbp]
	\centering
	\begin{minipage}[t]{0.48\linewidth}
		\centering
		\vspace{0pt}
		{\small\textbf{(a) Performance profiles}\par}
		\smallskip
		\includegraphics[width=\linewidth]{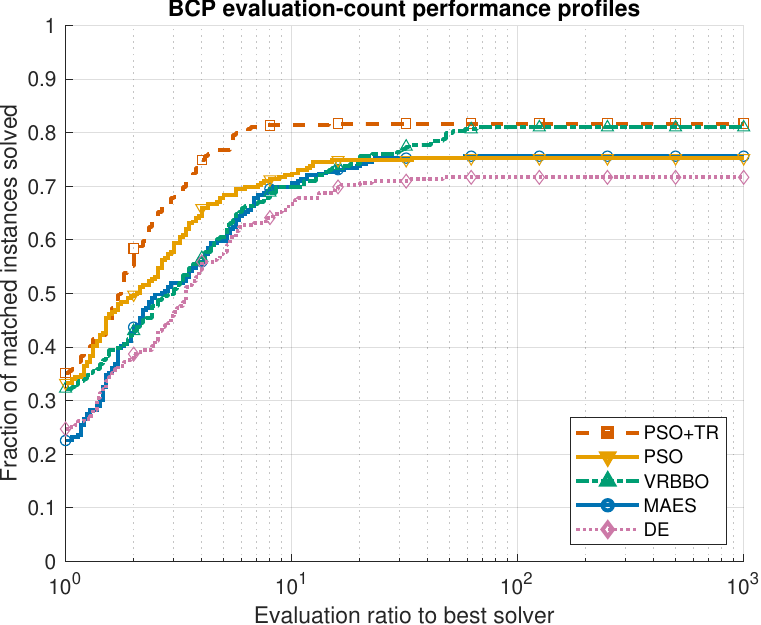}
	\end{minipage}
	\hfill
	\begin{minipage}[t]{0.48\linewidth}
		\centering
		\vspace{0pt}
		{\small\textbf{(b) Data profiles}\par}
		\smallskip
		\includegraphics[width=\linewidth]{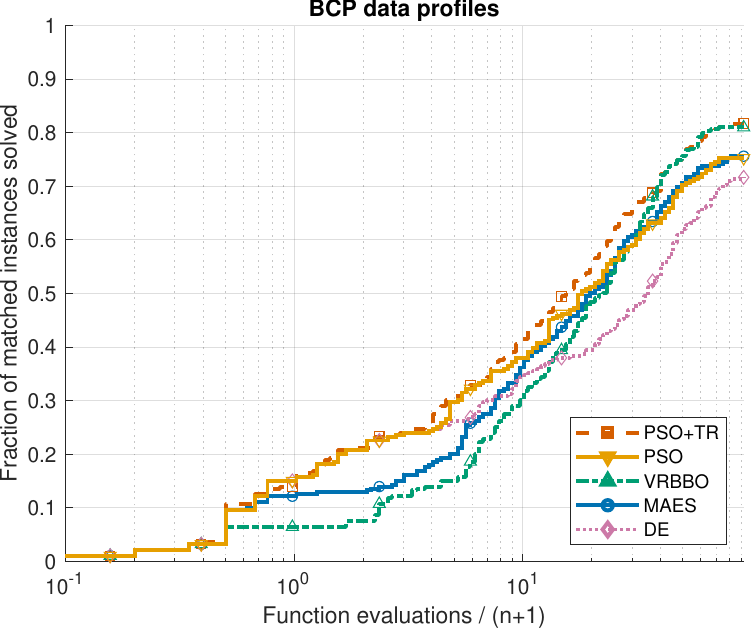}
	\end{minipage}
	\caption{Evaluation-count profiles for five DFO methods on
		279 matched \texttt{BCP} problem--seed instances
		(93 problems, three seeds).
		(a) Target attainment within a factor $\tau$ of the smallest
		observed first-hit count on each instance.
		(b) Target attainment within $\kappa(n_p+1)$ charged evaluations.
		Higher curves indicate greater attainment at the stated allowance.
		Unsolved instances remain in the common denominator, and all
		methods use matched evaluation budgets.}
	\label{fig:bcp-profiles}
\end{figure}

Figure~\ref{fig:bcp-profiles} compares relative first-hit efficiency
in panel~(a) with target attainment under a dimension-normalized
evaluation allowance in panel~(b).
At $\tau=1$, \texttt{PSO+TR} has the largest fraction of fastest
recorded first hits, including ties, at 35.13\%.
Both profiles approach the target-attainment fractions reported below.

\paragraph{Comparison results.}
Table~\ref{tab:bcp-five-summary} presents the five-solver results
on 279 matched problem--seed instances. All methods use identical
objective-evaluation limits on each matched instance, with targets
allowing 1\% of the initial-to-best observed objective gap.
The proposed \texttt{PSO+TR} reaches the common target on 81.72\%
of the instances, compared with 81.00\% for \texttt{VRBBO},
75.63\% for \texttt{MAES}, 75.27\% for \texttt{PSO}, and
71.68\% for \texttt{DE}.
Table~\ref{tab:bcp-psotr-vs-pso} summarizes the pairwise comparison,
with a target-attainment advantage of 6.45 percentage points
for \texttt{PSO+TR} over \texttt{PSO}.

The median first-hit evaluation counts, computed separately on
each method's successful instances, are 41.5 for \texttt{PSO+TR},
54 for \texttt{VRBBO}, 44 for \texttt{MAES}, 36 for \texttt{PSO},
and 46 for \texttt{DE}. Thus, the smaller conditional median
for \texttt{PSO} coexists with its lower overall
target-attainment fraction.

Mean regret averages instance-normalized regret values rather
than raw objective values. ``Tied best'' counts instances satisfying
$|f_{p,s}^{\rm final}-f_p^\star|\le
10^{-10}\max\{1,|f_p^\star|\}$; several methods may be counted
on the same instance. \texttt{VRBBO} records the largest number
of tied-best final values, whereas \texttt{PSO} has a smaller
mean regret than \texttt{PSO+TR}. The observed ranking therefore
depends on the criterion.

These differences are descriptive: three stochastic repetitions
of the same 93 problems do not constitute 279 independent problem
families, and no statistical significance or universal superiority
is inferred. Moreover, the separate \texttt{PSO} implementation
differs from \texttt{PSO+TR} in swarm initialization and
boundary-velocity handling, so their comparison does not isolate
the contribution of the \texttt{TR} safeguard.

\begin{table}[!htbp]\centering\small
	\caption{Five-solver \texttt{BCP} results on 279 matched
		problem--seed instances.}
	\label{tab:bcp-five-summary}
	\begin{tabular}{lrrrr}\toprule
		Method & Solved (\%) & Median evaluations & Mean regret & Tied best \\
		\midrule
		\texttt{MAES} & 75.63 & 44 & $1.494\times10^{4}$ & 59 \\
		\texttt{PSO+TR} & 81.72 & 41.5 & 15.27 & 89 \\
		\texttt{VRBBO} & 81.00 & 54 & $7.127\times10^{3}$ & 148 \\
		\texttt{DE} & 71.68 & 46 & 10.95 & 63 \\
		\texttt{PSO} & 75.27 & 36 & 0.9065 & 72 \\
		\bottomrule
	\end{tabular}
\end{table}

\begin{table}[!htbp]\centering\small
	\caption{\texttt{PSO+TR} and baseline \texttt{PSO} results
		on the same 279 matched \texttt{BCP} instances.}
	\label{tab:bcp-psotr-vs-pso}
	\begin{tabular}{lrrr}\toprule
		Method & Solved (\%) & Median evaluations & Tied best \\
		\midrule
		\texttt{PSO+TR} & 81.72 & 41.5 & 89 \\
		\texttt{PSO} & 75.27 & 36 & 72 \\
		\bottomrule
	\end{tabular}
\end{table}

\paragraph{Conditional \texttt{TR} and subspace diagnostics.}
Over the 279 matched instances, \texttt{PSO+TR} makes 1,409 conditional \texttt{TR}
invocations, accepts 1,687 \texttt{TR} steps, and expends 32,988 objective evaluations
on \texttt{TR} operations, including model construction. These quantities are
compatible with a successive-\texttt{TR} mechanism: one invocation may accept more
than one local step. Consequently, 119.7 accepted steps per 100 invocations
is a frequency, not a probability or percentage of accepted invocations.
The logs record 3,046 basis-column updates and 77 refreshes.
An update following an accepted \texttt{TR} step preserves the subspace
span and can only change the selected column's sign in exact arithmetic;
these update counts therefore do not count distinct changes of span.
The maximum observed orthogonality residual is
$\|S^\top S-I\|_F=6.99\times10^{-16}$. These data document the
safeguard's computational activity, but do not by themselves identify
its incremental benefit relative to an otherwise identical \texttt{PSO} run.

\begin{table}[htbp]\centering\small
	\caption{\texttt{PSO+TR} safeguard diagnostics for the \texttt{BCP} comparison.}\label{tab:bcp-psotr-diagnostics}
	\begin{tabular}{lr}\toprule
		Metric & Value \\ 
		\midrule
		\texttt{TR} invocations & 1409 \\ 
		Accepted \texttt{TR} steps & 1687 \\ 
		Accepted \texttt{TR} steps per 100 invocations & 119.7 \\ 
		\texttt{TR} objective evaluations & 32988 \\ 
		Subspace updates & 3046 \\ 
		Subspace refreshes & 77 \\ 
		Maximum $\|S^\top S-I\|_F$ & 6.99e-16 \\ 
		\bottomrule\end{tabular}\end{table}

\begin{figure}[!htbp]\centering
	\includegraphics[width=.87\linewidth]{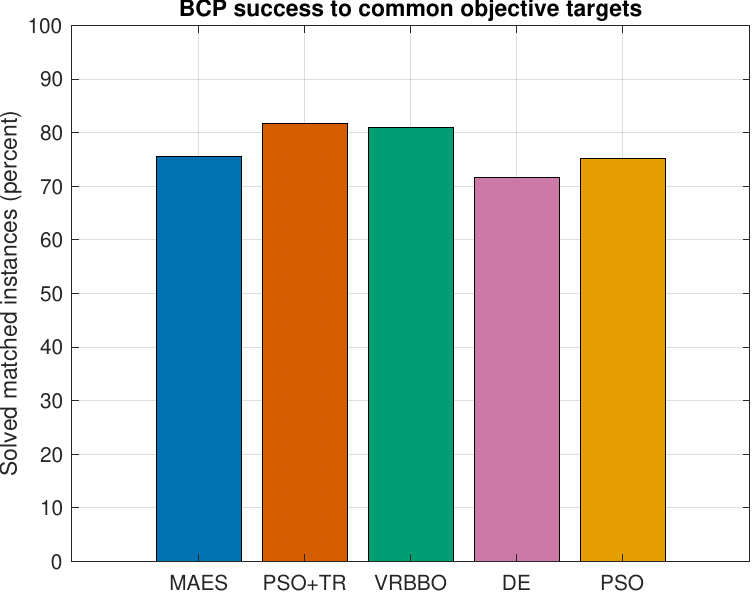}
	\caption{Observed proportion of the 279 matched \texttt{BCP} instances meeting the
		same instance-specific objective target for each of the five methods.}
	\label{fig:bcp-target-success}\end{figure}

\begin{figure}[!htbp]\centering
	\includegraphics[width=.87\linewidth]{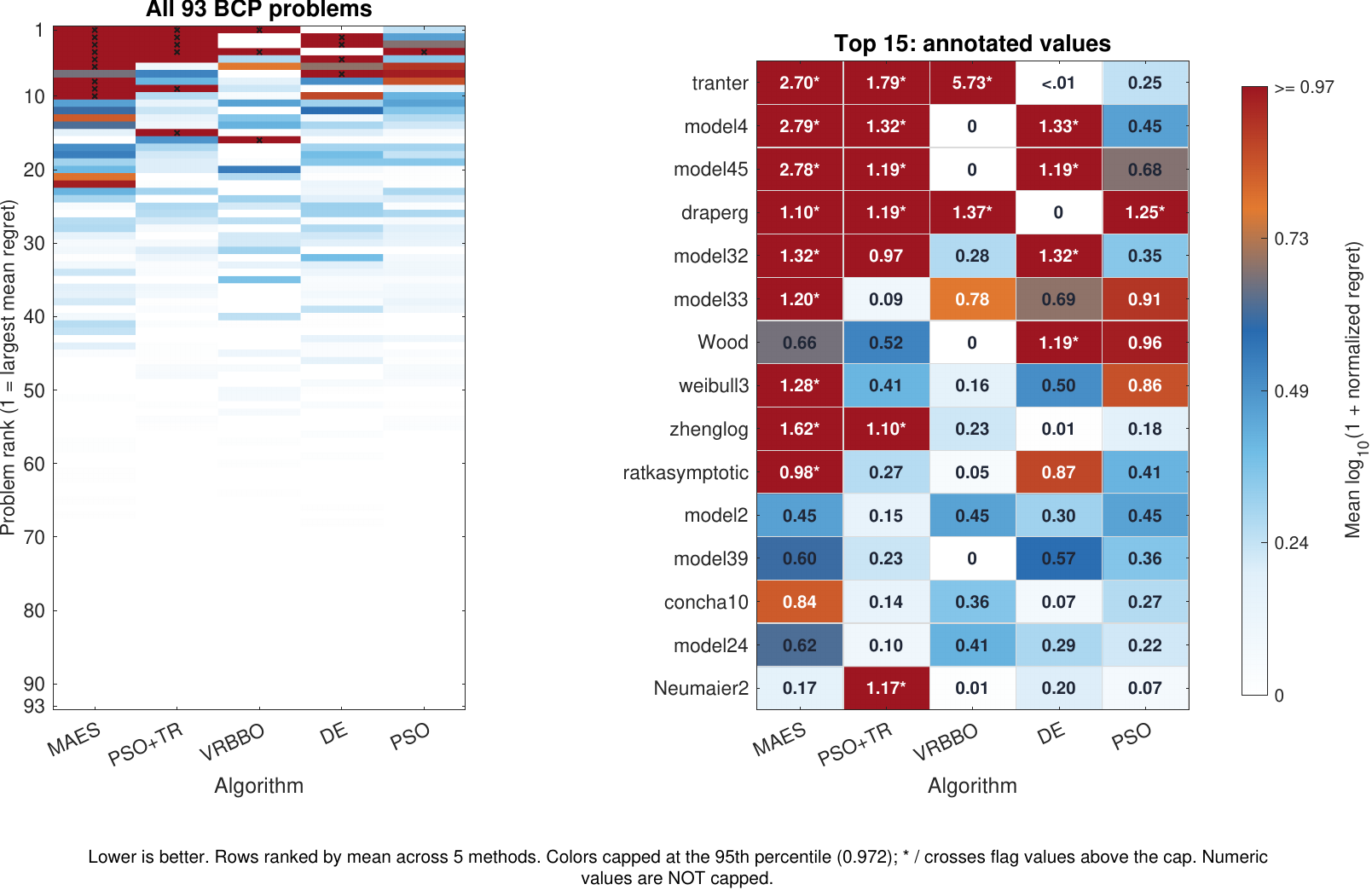}
	\caption{Problem-level mean $\log_{10}(1+\text{normalized regret})$ over
		three seeds. Each row is one of 93 eligible \texttt{BCP} problems; \texttt{MAES} and \texttt{PSO+TR}
		remain separate methods.}\label{fig:bcp-regret-heatmap}\end{figure}

Table~\ref{tab:bcp-psotr-diagnostics} collects these intervention
counts and numerical orthogonality diagnostics. Figure~\ref{fig:bcp-target-success} displays the target-attainment
fractions: \texttt{PSO+TR} solves 228 instances and \texttt{VRBBO}
solves 226, so their observed difference corresponds to two instances.
Figure~\ref{fig:bcp-regret-heatmap} exposes problem-level variation
in the mean transformed regret; the lowest plotted values occur for
different methods on different problems.

\paragraph{Exploratory hyperparameter sensitivity.}
The accompanying sensitivity figures summarize a separate exploratory
training grid in which the initial \texttt{TR} radius and subspace refresh
interval were varied, while the other displayed parameters were held fixed.
The saved parameter-selection record reports a budget factor of
1,000, whereas the final \texttt{BCP} comparison uses a factor of 100.
The plotted grid illustrates sensitivity within its own experimental
setting; its fixed-parameter panels provide no information about varying
those parameters. No final-comparison seed was used to select the
configuration reported above.

\begin{figure}[!htbp]\centering
	\includegraphics[width=.87\linewidth]{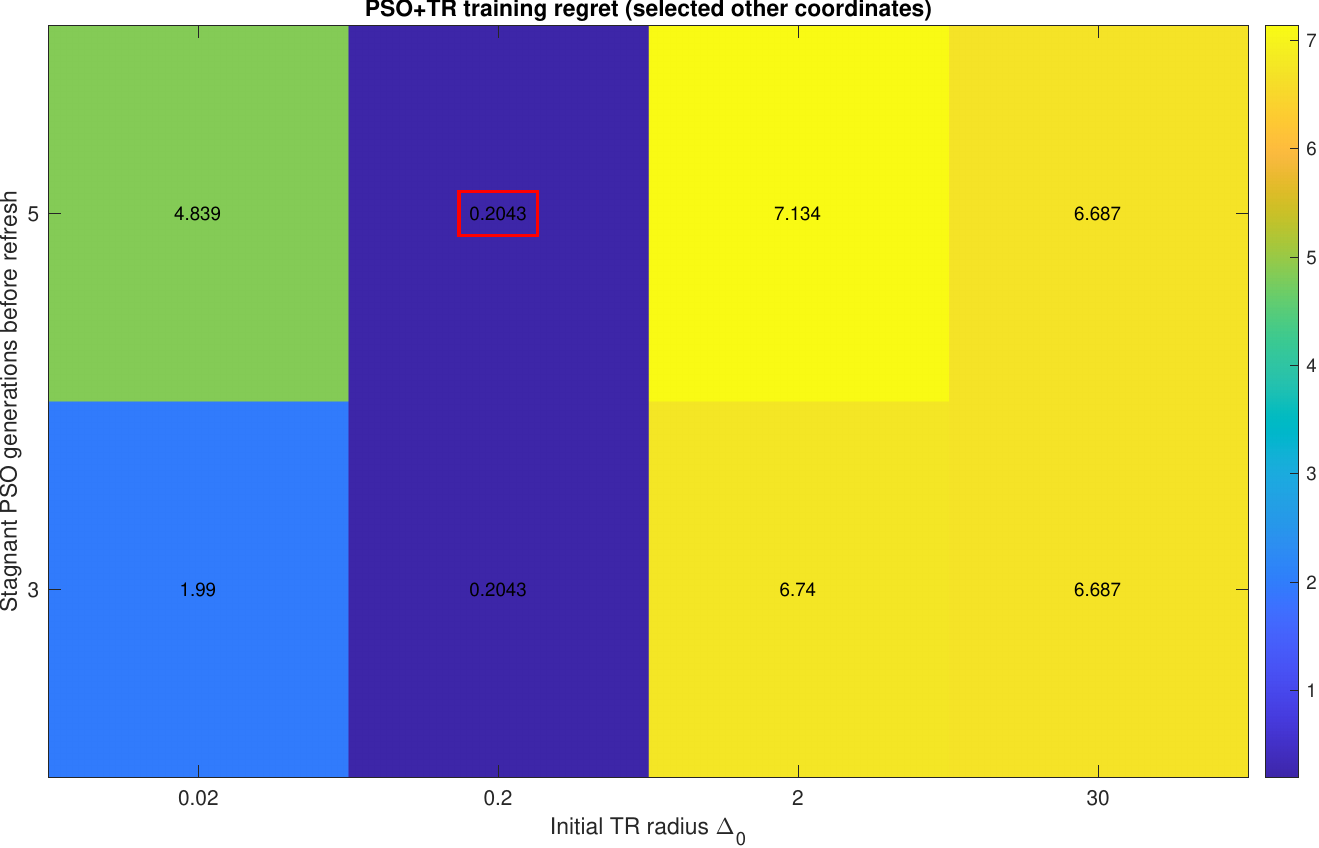}
	\caption{\texttt{PSO+TR} training normalized-regret grid as a function of the
		initial \texttt{TR} radius and subspace refresh interval; other grid coordinates
		are held at their selected values.}\label{fig:bcp-tuning-heatmap}\end{figure}

\begin{figure}[!htbp]\centering
	\includegraphics[width=.87\linewidth]{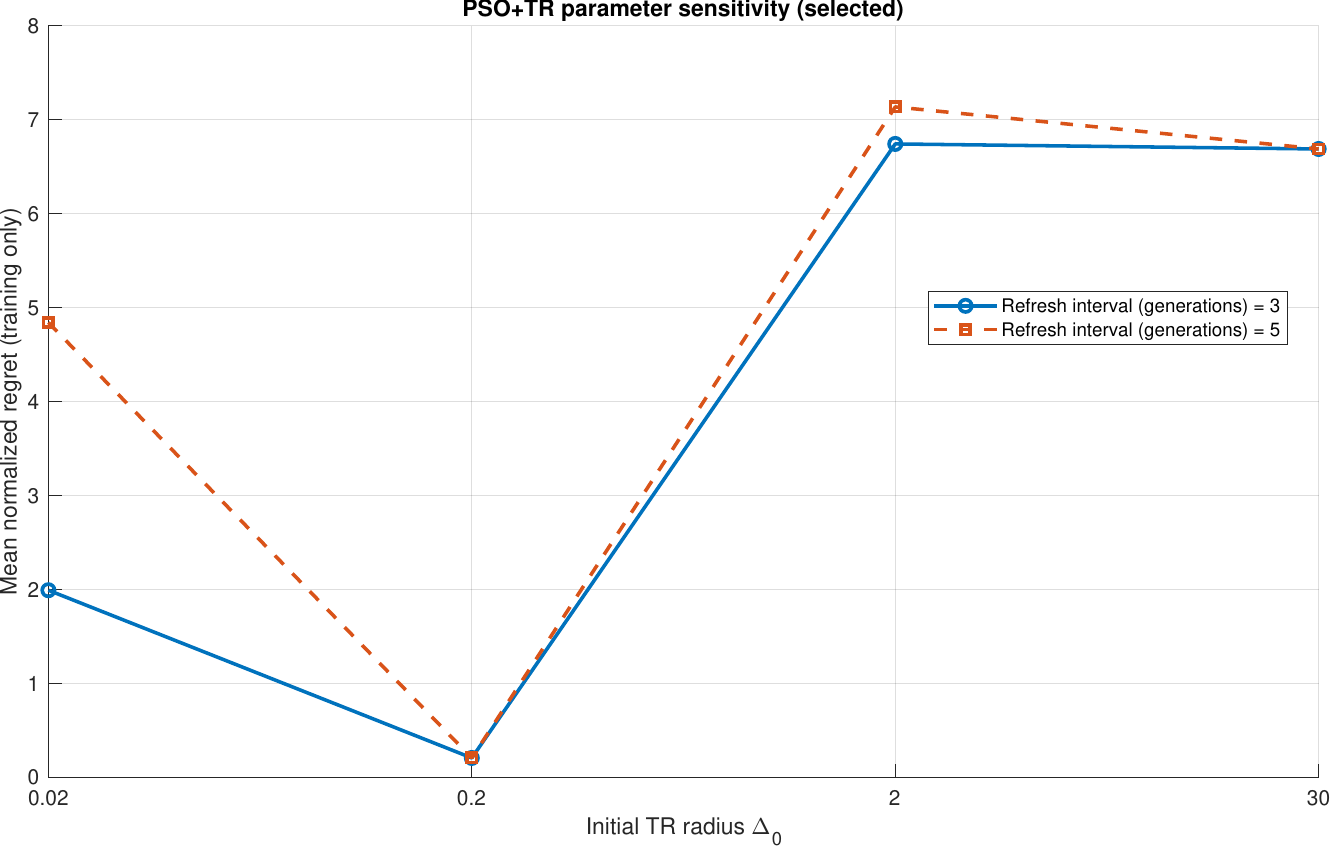}
	\caption{\texttt{PSO+TR} training-grid sensitivity, without holdout-based model
		selection.}\label{fig:bcp-tuning-sensitivity}\end{figure}

\begin{figure}[!htbp]\centering
	\includegraphics[width=.93\linewidth]{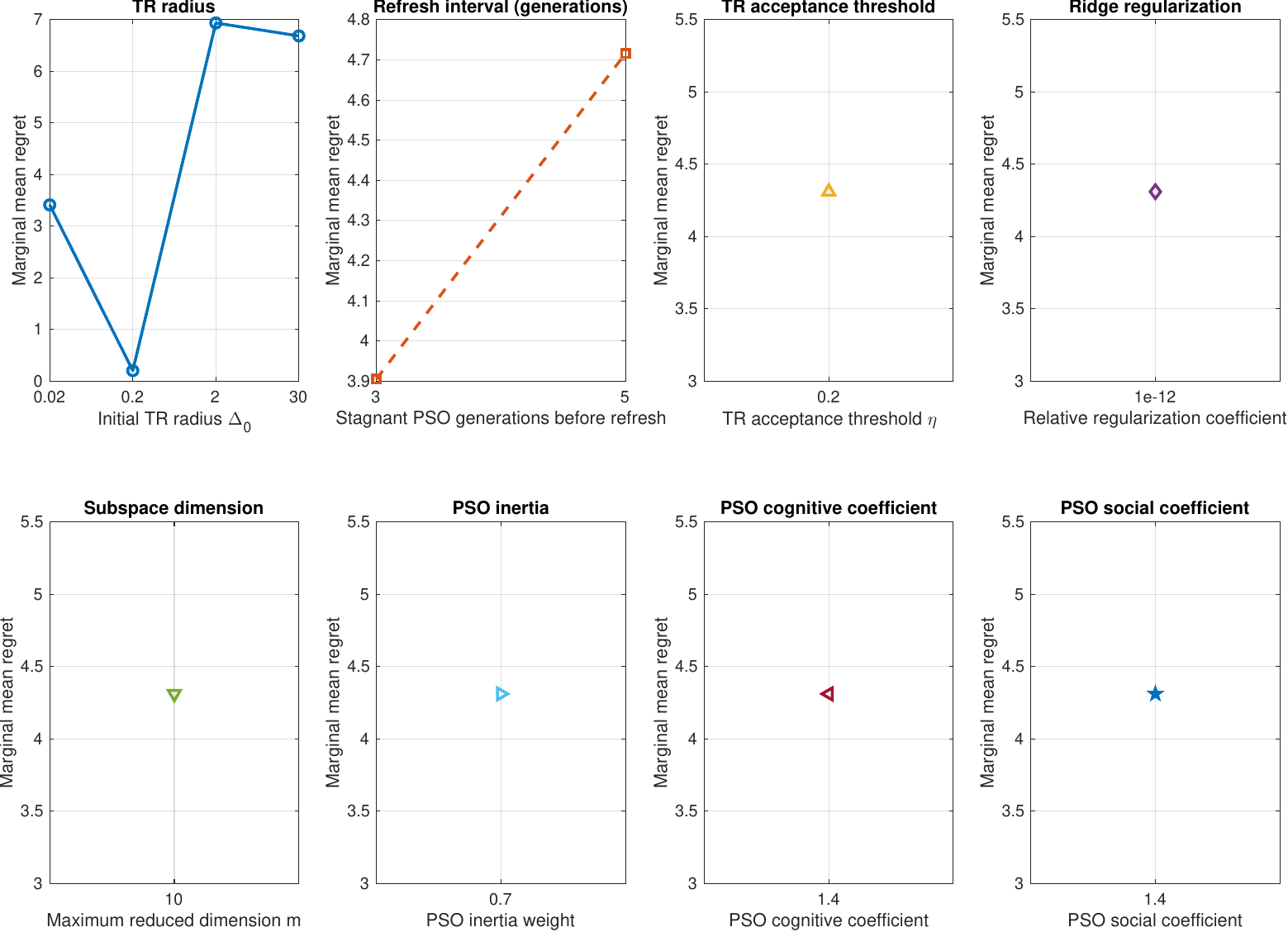}
	\caption{Marginal training-grid regret by parameter. The six parameters
		held fixed in the archived eight-configuration grid have only one value,
		so those panels are not evidence of sensitivity to those parameters.}
	\label{fig:bcp-tuning-marginal}\end{figure}

Figures~\ref{fig:bcp-tuning-heatmap} and~\ref{fig:bcp-tuning-sensitivity}
show the lowest mean training regret at the tested initial radius
$\Delta_0=0.2$ for both tested refresh intervals, three and five
generations. Figure~\ref{fig:bcp-tuning-marginal} presents the
coordinate marginals; its six fixed-parameter panels describe the
tested settings and do not establish sensitivity to those parameters.

\FloatBarrier
\section{Discussion and Conclusion}
\label{sec:discussion}

To distinguish supported conclusions from remaining questions,
we interpret the \texttt{BCP} optimizer comparison and robotic
controller results together and identify the controlled ablations
still needed. The benchmark evaluates target attainment and
evaluation efficiency, while the robotic experiment assesses
the complete tuning and deployment pipeline.

\texttt{PSO+TR} combines population exploration with conditional,
successive reduced-subspace \texttt{TR} steps. Its persistent basis and
regularized swarm-geometry surrogate support local refinement
without constructing a full-dimensional curvature matrix.
All attempted evaluations are charged to a common objective-evaluation
budget. Controlled ablations are needed to isolate the contributions
of the \texttt{TR} safeguard, subspace updates, and curvature construction.

On the 93 \texttt{BCP} problems, \texttt{PSO+TR} achieves the
highest target-attainment rate among the five compared
methods, reaching the prescribed target on 81.72\% of the
279 matched instances, compared with 75.27\% for
\texttt{PSO} and 81.00\% for \texttt{VRBBO}. Nevertheless,
\texttt{PSO} has a lower median first-hit evaluation count
among successful instances and a lower mean normalized
regret, indicating that no method dominates across all
performance measures. The \texttt{TR} safeguard records 1,409
invocations and 1,687 accepted steps, consuming 32,988
objective evaluations. These findings demonstrate the
practical occurrence of successive local improvements,
although differences between the compared implementations
prevent attributing the overall performance gain solely
to the \texttt{TR} component.

In the robotic DFPO study, the validated \texttt{PSO+TR}
deployment pipeline achieves completion rates of 97.78\%
and 58.00\% under development and withheld stress
conditions, respectively, compared with 91.22\% and
31.12\% for the fixed staged controller. Under stress
conditions, the violation rate decreases from 20.69\%
to 2.19\%, while the mean episode cost decreases from
13.124 to 10.094. Seven of eight optimized controllers
pass independent validation. However, only two local
steps are accepted across 11 robotic \texttt{TR} invocations.
Consequently, the observed improvements characterize
the complete tuning and validation pipeline rather than
the isolated contribution of the \texttt{TR} safeguard.

Overall, the results establish competitive performance
on the considered DFO benchmarks and demonstrate the
potential of \texttt{PSO+TR} for tuning interpretable
robotic insertion controllers. Nevertheless, the robotic
experiments rely on a synthetic planar contact model,
and their results do not establish physical-robot
performance or safety guarantees. Future research will
focus on controlled \texttt{PSO}-only ablations,
sensitivity to algorithmic parameters and contact
uncertainty, more realistic three-dimensional assembly
models, and experimental validation on physical robots.

\appendix
\section{Seed-level outcomes of the robotic controller-tuning experiment}
\label{app:psotr-seed-results}

To expose variation hidden by aggregate results,
Table~\ref{tab:psotr-native-seeds} reports completion and violation
counts, mean episode cost, and mean termination time for each
of the eight independent optimizer seeds. Within each seed,
the four controller configurations use matched final-test scenarios,
with 400 episodes per distribution and controller.

Seed~4 produced no eligible screening candidate and therefore
deployed the fixed staged reference without candidate validation.
Its reported deployed-policy outcomes are the previously observed
fixed staged outcomes on those same paired scenarios.
No additional simulation episodes were generated for this correction.

Averaging the deployed-policy completion fractions across the eight
seeds gives 97.78\% for the development distribution and 58.00\%
for the stress distribution, consistent with the aggregate results.
Variation across seeds combines differences in the independent
optimization runs and in the separately drawn final-test scenarios;
it therefore does not isolate variability due to optimization alone.

\begin{longtable}{llrrrr}
	\caption{Native \texttt{PSO+TR} deployed-policy results for all eight seeds on the reconstructed simulator (400 paired cases per group and seed). Here, mean time is the average
		recorded termination time $T_e$ over all 400 episodes, including
		timeouts and protective terminations.}\label{tab:psotr-native-seeds}\\
	\toprule Group & Seed & Successes & Violations & Mean cost & Mean time (s) \\ 
	\midrule
	\endfirsthead
	\toprule Group & Seed & Successes & Violations & Mean cost & Mean time (s) \\ 
	\midrule
	\endhead
	nominal & 0 & 400 & 0 & 2.709 & 2.165 \\ 
	nominal & 1 & 400 & 0 & 3.263 & 2.585 \\ 
	nominal & 2 & 400 & 0 & 2.835 & 2.447 \\ 
	nominal & 3 & 400 & 0 & 2.930 & 2.373 \\ 
	nominal & 4 & 400 & 0 & 4.401 & 3.547 \\ 
	nominal & 5 & 400 & 0 & 2.573 & 2.105 \\ 
	nominal & 6 & 400 & 0 & 2.795 & 2.268 \\ 
	nominal & 7 & 400 & 0 & 2.890 & 2.471 \\ 
	development & 0 & 398 & 0 & 3.591 & 3.290 \\ 
	development & 1 & 396 & 0 & 4.038 & 3.334 \\ 
	development & 2 & 400 & 0 & 3.573 & 3.589 \\ 
	development & 3 & 398 & 0 & 3.952 & 3.795 \\ 
	development & 4 & 358 & 0 & 6.650 & 5.351 \\ 
	development & 5 & 394 & 0 & 3.613 & 3.409 \\ 
	development & 6 & 386 & 0 & 4.019 & 3.409 \\ 
	development & 7 & 399 & 0 & 3.686 & 3.638 \\ 
	stress & 0 & 233 & 0 & 10.333 & 5.914 \\ 
	stress & 1 & 197 & 2 & 11.836 & 6.412 \\ 
	stress & 2 & 265 & 0 & 8.634 & 5.439 \\ 
	stress & 3 & 253 & 0 & 9.346 & 5.639 \\ 
	stress & 4 & 119 & 68 & 13.467 & 7.226 \\ 
	stress & 5 & 282 & 0 & 7.988 & 4.858 \\ 
	stress & 6 & 246 & 0 & 9.946 & 5.753 \\ 
	stress & 7 & 261 & 0 & 9.202 & 5.732 \\ 
	\bottomrule\end{longtable}

\end{sloppypar}

\section*{Funding}

M.K. acknowledges financial support of the Austrian Science Fund (FWF) under \url{https://doi.org/10.55776/PAT2747625}.


\end{document}